\documentclass[12pt,a4paper]{article}

\usepackage{amsmath, amssymb, amsthm}
\usepackage[dvips,final]{epsfig}
\usepackage[dvips,final]{graphicx}
\usepackage{enumitem}

\usepackage{color}

\definecolor{green1}{rgb}{0.1,0.65,0}
\definecolor{green2}{rgb}{0.2,0.5,0}
\definecolor{viol}{rgb}{0.5,0,0.5}
\definecolor{cyan1}{rgb}{0,0.7,0.3}

\usepackage{titlesec}
\titleformat{\section}{\bfseries}{\thesection}{1em}{}
\titleformat{\subsection}{\itshape}{\thesubsection}{1em}{}

\newfont{\ctv}{msam10}

\def\ve{\varepsilon}

\def\dd{\,\mathrm{d}}
\def\dive{\mathrm{\,div\,}}

\def\supess{\mathop{\mathrm{sup\,ess}}}

\def\ani{\mbox{\ctv \symbol{75}}}
\def\dch{(\chi{+}\rho^*(1{-}\chi))}
\def\water{\mathrm{H}_2\mathrm{O}}
\def\e{\mathrm{e}}

\def\real{\mathbb{R}}
\def\nat{\mathbb{N}}
\def\tens{\real^{3\times 3}_{\mathrm{sym}}}

\def\AA{\mathbf{A}}
\def\BB{\mathbf{B}}

\def\CC{\mathcal{C}}

\def\Ucal{\mathcal{U}}
\def\Wcal{\mathcal{W}}
\def\play{\mathfrak{f}}

\def\nas{\nabla_s}

\def\io{\int_{\Omega}}
\def\ipo{\int_{\partial\Omega}}

\def\softd{{\leavevmode\setbox1=\hbox{d}%
\hbox to 1.05\wd1{d\kern-0.4ex{\char039}\hss}}}

\newcommand{\for}{\ \mbox{for } \,}

\def\be{\begin{equation}\label}
\def\ee{\end{equation}}

\newcommand{\bbox}{\mbox{\ctv \symbol{4}}}
\def\QED{{${}\hfill\bbox$}}
\newenvironment{pf}[1]{\par\vskip1mm{\noindent\it #1.}\ }{\QED\par
\vskip2mm}

\def\bpf{\begin{pf}}
\def\epf{\end{pf}}

\newtheorem{theorem}{Theorem}[section]

\newtheorem{hypo}[theorem]{Hypothesis}
\newtheorem{propo}[theorem]{Proposition}

\newtheorem{remark}[theorem]{Remark}

\newenvironment{pkrev}{\color{blue}}{\color{black}}
\newcommand{\bpk}{\begin{pkrev}}
\newcommand{\epk}{\end{pkrev}}

\numberwithin{equation}{section}

\begin{document}

\title{Phase transitions in porous media\thanks{The support from the Austrian Science Fund (FWF) projects V662 and F65, from the GA\v CR Grant No.~20-14736S, and from the European Regional Development Fund, Project No. CZ.02.1.01/0.0/0.0/16{\_}019/0000778 is gratefully acknowledged.}}

\author{Chiara Gavioli
\thanks{Institute of Analysis and Scientific Computing, TU Wien, Wiedner Hauptstra\ss e 8-10, A-1040 Vienna, Austria, E-mail: {\tt chiara.gavioli@tuwien.ac.at}.}
\and Pavel Krej\v c\'{\i}
\thanks{Faculty of Civil Engineering, Czech Technical University, Th\'akurova 7, CZ-16629 Praha 6, Czech Republic, E-mail: {\tt Pavel.Krejci@cvut.cz}.}
}

\date{}

\maketitle

\begin{abstract}
The full quasistatic thermomechanical system of PDEs, describing water diffusion with the possibility of freezing and melting in a visco-elasto-plastic porous solid, is studied in detail under the hypothesis that the pressure-saturation hysteresis relation is given in terms of the Preisach hysteresis operator. The resulting system of balance equations for mass, momentum, and energy coupled with the phase dynamics equation is shown to admit a global solution under general assumptions on the data.

\bigskip


\medskip


\end{abstract}


\section*{Introduction}\label{int}

A model for fluid flow in partially saturated porous media with thermomechanical interaction
was proposed and analyzed in \cite{ak,dkr}.
The model was subsequently extended in \cite{krs} by including the effects of freezing and melting of the water in the pores. Typical examples in which such situations arise are related to groundwater flows and to the freezing-melting cycles of water sucked into the pores of concrete. Due to the specific volume difference between water and ice, this process produces important pressure changes and represents one of the main reasons for the degradation of construction materials in buildings, bridges, and roads. The model of \cite{krs} still neglects the influence of changes of microstructure, as for example the breaking of pores, but the main thermomechanical interactions between the state variables are taken into account.    

The model is based on the assumption that slow diffusion of the fluid through the porous solid is a dominant effect, so that the Lagrangian description is considered to be appropriate. It is assumed that volume changes of the solid matrix material are negligible with respect to the pore volume evolution during the process. The pores are filled with a mixture of $\water$ and gas, and $\water$ itself is a mixture of the liquid (water) and the solid phase (ice). That is, in addition to the standard state variables like capillary pressure, displacement, and absolute temperature, we need to consider the evolution of a phase parameter $\chi$ representing the relative proportion of water in the $\water$ part and its influence on pressure changes due to the different mass densities of water and ice.

The resulting system consists of mechanical balance equation for the deformations of the solid body, mass balance equation based on the Darcy law for the fluid diffusion with interaction terms similar to the Biot system studied, e.\,g., in \cite{show}, a differential inclusion for the phase fraction $\chi$ of relaxed Stefan type as in \cite{vis} and governing the water-ice phase transition, and the energy balance equation derived from the first and the second principles of thermodynamics with heat sources due to viscosity, plasticity, diffusion, and phase transition. 

The present paper develops the ideas of \cite{krs} in the sense that the effects of capillary hysteresis, which is assumed to be of Preisach type in agreement with the results of \cite{flynn}, shear stresses, and elastoplasticity are considered in full generality. This represents an enormous increase of mathematical complexity. While the momentum balance in the shear-stress-free case in \cite{krs} can be reduced to an ODE, here, we need to exploit deeper results from the theory of PDEs to control the interactions between individual components of the system, as well as a generalized Moser iteration scheme from \cite{gk}. Additional difficulties are due to the effects of the three heat sources produced by mechanical hysteresis dissipation (plasticity, capillarity, phase transitions).

The paper is divided into five sections. In the next Section~\ref{mod}, we briefly recall the principles of the model introduced in \cite{krs}, taking into account capillary hysteresis and elastoplastic hysteresis effects as in \cite{dkr}. Section \ref{sec_prei} is devoted to a short survey of the Preisach hysteresis model. In Section~\ref{sta}, we state the mathematical problem, the main assumptions on the data, and the main Theorem \ref{t1}, the proof of which is split into Sections~\ref{cutoff} and \ref{estiR}. The steps of the proof are as follows. In Section~\ref{cutoff}, we first cut off some of the pressure and temperature dependent terms in the system by means of a cut-off parameter $R$, regularize the mass balance equation with a fourth order term depending on an additional small regularizing parameter $\eta$, solve the related problem employing a Galerkin approximation scheme, and pass to the singular limit $\eta \to 0$ in the regularizing term. Then, in Section~\ref{estiR}, we derive a series of $R$-independent estimates like the energy estimate, the so-called Dafermos estimate (with negative small powers of the temperature), Moser-type and then higher-order estimates for the capillary pressure and for the temperature which allow us in Section \ref{sec_lim} to pass to the limit in the cut-off system as $R\to\infty$, which will conclude the proof of the existence result.


\section{The model}\label{mod}

Consider a bounded domain $\Omega\subset \real^3$ of class $C^{1,1}$ filled with an elastoplastic solid matrix material with pores containing a mixture of $\water$ and gas, where we assume that $\water$ may appear in one of the two phases: water or ice. We state the balance laws in referential (Lagrangian) coordinates. We have in mind construction materials where large deformations are not expected to occur. This hypothesis enables us to reduce the complexity of the problem and assume that the deformations are small in order to avoid higher degree nonlinearities. We denote for $x\in \Omega$ and time $t\in [0,T]$
\begin{description}
\item $p(x,t)$ ... capillary pressure;
\item $u(x,t)$ ... displacement vector in the solid;
\item $\ve(x,t) = \nas u(x,t)$ ... linear strain tensor, $(\nas u)_{ij} :=
\frac12 \left(\frac{\partial u_i}{\partial x_j} + \frac{\partial u_j}{\partial x_i}\right)$;
\item $\theta(x,t)$ ... absolute temperature;
\item $\chi(x,t) \in [0,1]$ ... relative amount of water in the $\water$ part.
\end{description}

The model derived in \cite{krs} aims at coupling the effects of capillarity, interaction between a deformable solid matrix material and $\water$ in the pores which may undergo water-ice phase transitions, and energy exchange between the individual components of the system. Hysteresis is included following the modeling section of \cite{dkr}. The full system consists of equations describing mass balance \eqref{e1}, mechanical equilibrium \eqref{e2}, energy balance \eqref{e3} and phase evolution \eqref{e4} in the form
\begin{align}
\big(\dch(G[p] + \dive u)\big)_t &= \dive (\mu(p)\nabla p), \label{e1} \\
-\dive(\BB\nas u_t + P[\nas u]) +\nabla(p \dch + \beta(\theta-\theta_c)) &= g, \label{e2} \\
\CC_V(\theta)_t - \dive(\kappa(\theta)\nabla\theta) = \BB \nas u_t:\nas u_t + \mu(p)|\nabla p|^2 &+ \|D_P[\nas u]_t\|_* \nonumber\\
+ \dch|D_G[p]_t| + \gamma(\theta,\dive u)\chi_t^2 &- \frac{L}{\theta_c}\theta\chi_t - \beta \theta \dive u_t, \label{e3} \\
-\gamma(\theta,\dive u) \chi_t + 
(1-\rho^*)\left(p\,G[p] - U_G[p] + p\dive u\right) + L\left(\frac{\theta}{\theta_c}-1\right) &\in \partial I_{[0,1]}(\chi). \label{e4}
\end{align}
We refer to \cite{krs} for the details of the physical arguments. Let us just mention that
the mass balance \eqref{e1} is derived from Darcy's law, and $\mu(p)\nabla p$ is the specific liquid mass flux. The constant $\rho^*\in (0,1)$ is the ratio between ice and water mass densities, whereas the symbol $G$ describes the pressure-saturation curve and, following \cite{ak,dkr}, is of the form
$$G[p] = f(p)+G_0[p].$$
Here $f$ is a bounded monotone function satisfying Hypothesis \ref{h1}\,(vi) below, whereas $G_0$ is the Preisach hysteresis operator from Section \ref{sec_prei}. The momentum balance \eqref{e2} is derived by assuming that the process is quasi-static, so that the inertia term is negligible. Here $\BB$ is a positive definite viscosity matrix, $P$ the constitutive operator of elastoplasticity defined below in \eqref{plast}, $\beta$ is the thermal expansion coefficient, $\theta_c$ is the melting temperature at standard pressure and $g$ is a given volume force (gravity, e.\,g.). The term $p \dch$ represents the pressure component due to the phase transition. Finally, the energy balance \eqref{e3} and the inclusion \eqref{e4} for the evolution of the phase parameter are derived from the principles of thermodynamics with the aid of the energy balance for both the plasticity and the pressure-saturation operator
\begin{equation}\label{en_idG}
	P[\ve]:\ve_t = U_P[\ve]_t + \|D_P[\ve]_t\|_*, \qquad	G[p]_t p = U_G[p]_t+|D_G[p]_t|,
\end{equation}
where $U_P,U_G$ and $D_P,D_G$ are the potential and dissipation operators, and $\|\cdot\|_*$ is a seminorm in the space $\tens$ of symmetric $3\times 3$ tensors. In \eqref{e3}, $\CC_V(\theta)$ is the caloric component of the internal energy, $\kappa(\theta)$ is the heat conductivity coefficient, $L$ is the latent heat, $I_{[0,1]}$ is the indicator function of the interval $[0,1]$ and $\partial I_{[0,1]}$ is its subdifferential, $\gamma(\theta,\dive u)$ is the phase relaxation time (which we assume to explicitly depend on both $\theta$ and $\dive u$ for technical reasons).

Note that the values of $G$ have to be naturally confined between $0$ and $1$, so that the system is degenerate in the sense that we do not control a priori the time derivatives of $p$ in \eqref{e1}. Another difficulty is related to the lack of spatial regularity of $\chi$. The temperature field is problematic as well: equation \eqref{e3} contains high order heat source terms, which are difficult to handle and prevent the temperature from being regular.

We complement the system \eqref{e1}--\eqref{e4} with initial conditions
\begin{equation}\label{ini}
	\left.
	\begin{array}{rcl}
		p(x,0) &=& p^0(x) \\
		u(x,0) &=& u^0(x) \\
		\theta(x,0) &=& \theta^0(x) \\
		\chi(x,0) &=& \chi^0(x)
	\end{array}
	\right\rbrace \mbox{in}\ \Omega\,,
\end{equation}
and boundary conditions
\begin{equation}\label{boufull}
	\left.
	\begin{array}{rcl}
		u &=& 0 \\
		\mu(p)\nabla p\cdot n &=& \alpha(x)(p^*-p) \\
		\kappa(\theta)\nabla\theta\cdot n &=& \omega(x)(\theta^*-\theta)
	\end{array}
	\right\rbrace \mbox{on}\ \partial\Omega\,,
\end{equation}
where $p^*$ is a given outer pressure, $\theta^*$ is a given outer temperature, $\alpha : \partial\Omega \to [0,\infty)$ is the permeability of the boundary and $\omega : \partial\Omega \to [0,\infty)$ is the heat conductivity of the boundary.

The solution to \eqref{e1}--\eqref{e4} was only constructed in \cite{krs} under the assumption that shear stresses in the momentum balance equation \eqref{e2} as well as all hysteresis effects are neglected. Then \eqref{e2} turns into an ODE for the relative volume change $\dive u$, which considerably simplifies the analysis. Here we prove existence of a global solution for the full problem under suitable hypotheses.


\section{Hysteresis in capillarity phenomena}\label{sec_prei}

The operator $G$ is considered as a sum
\begin{equation}\label{G}
	G[p]=f(p)+G_0[p],	
\end{equation}
where $f$ is a monotone function satisfying Hypothesis \ref{h1}\,(vi) in Section \ref{sta} below, and $G_0$ is a Preisach operator that we briefly describe here.

For a given input function $p\in W^{1,1}(0,T)$ and a memory parameter $r>0$ we define the scalar function $\xi_r(t)$ as the solution of the variational inequality
\begin{equation}
	\begin{cases}
		|p(t)-\xi_r(t)|\le r &\forall t\in[0,T],\\
		(\xi_r)_t(p(t)-\xi_r(t)-z)\ge 0 \ &\text{a.\,e. }\forall z\in[-r,r],
	\end{cases}
	\label{play}
\end{equation}
with prescribed initial condition 
\begin{equation}\label{ini_play}
	\xi_r(0)=\max\{p(0)-r,\min\{0,p(0)+r\}\}.
\end{equation}
We have indeed for all $r>0$ the initial bound
\begin{equation}\label{ini_bound}
	\xi_r(0)\le|p(0)|.
\end{equation}
The mapping $\play_r:W^{1,1}(0,T)\rightarrow W^{1,1}(0,T)$ which with each $p\in W^{1,1}(0,T)$ associates the solution $\xi_r=\play_r[p]\in W^{1,1}(0,T)$ of \eqref{play}--\eqref{ini_play} is called the \textit{play}. This concept goes back to \cite{krpo}, and the proof of the following statements can be found e.g. in \cite{kb}.
\begin{propo}
	For each $r>0$, the mapping $\play_r:W^{1,1}(0,T)\rightarrow W^{1,1}(0,T)$ is Lipschitz continuous and admits a Lipschitz continuous extension to $\play_r:C[0,T]\rightarrow C[0,T]$ in the sense that for every $p_1,p_2\in C[0,T]$ and every $t\in[0,T]$ we have
	\begin{equation}\label{lip_frak}
		|\play_r[p_1](t)-\play_r[p_2](t)|\le\max_{\tau\in[0,t]}|p_1(\tau)-p_2(\tau)|.
	\end{equation}
	Moreover, for each $p\in W^{1,1}(0,T)$, the energy balance equation
	\begin{equation}\label{ene_frak}
		\play_r[p]_tp-\frac{1}{2}\left(\play_r^2[p]\right)_t=|r\play_r[p]_t|
	\end{equation}
	and the identity
	\begin{equation}\label{id_frak}
		\play_r[p]_tp_t=\left(\play_r[p]_t\right)^2
	\end{equation}
	hold a.\,e. in $(0,T)$.
\end{propo}

Given a nonnegative function $\psi\in L^1((0,\infty)\times\real)$, called the \textit{Preisach density}, we define the Preisach operator $G_0$ as a mapping that with each $p\in C[0,T]$ associates the integral
\begin{equation}\label{G0}
	G_0[p](t)=\int_0^\infty\int_0^{\play_r[p](t)}\psi(r,v)\dd v\dd r.
\end{equation}
Hence,
\be{G0a}
G_0[p](t) \le \int_0^\infty\int_{0}^\infty \psi(r,v)\dd v\dd r =: C_\psi^+, \quad - G_0[p](t) \le \int_0^\infty\int_{-\infty}^0 \psi(r,v)\dd v\dd r =: C_\psi^-
\ee
for all $p \in C[0,T]$ and all $t \in [0,T]$, and we assume
\be{G0b}
0 < C_\psi^\pm < \frac12.
\ee
From \eqref{ene_frak}--\eqref{G0} we immediately deduce the Preisach energy identity
\begin{equation}\label{en_id0}
	G_0[p]_tp-U_0[p]_t=|D_0[p]_t|\text{ a.\,e.},
\end{equation}
provided we define the Preisach potential $U_0$ and the dissipation operator $D_0$ by the integrals
\begin{equation}
	U_0[p](t)=\int_0^\infty\int_0^{\play_r[p](t)}v\psi(r,v)\dd v\dd r, \qquad D_0[p](t)=\int_0^\infty\int_0^{\play_r[p](t)}r\psi(r,v)\dd v\dd r.
	\label{U0_D0}
\end{equation}
The energy identity in \eqref{en_idG} then holds with the choice
\begin{equation}\label{U_D}
	U_G[p] = pf(p)-\int_0^p f(z)\dd z+U_0[p], \qquad D_G[p] = D_0[p].
\end{equation}
With the notation
\begin{equation}\label{Phi_V}
	\Phi(p)=\int_0^p f(z)\dd z, \quad V(p)=pf(p)-\Phi(p) = \int_0^p f'(z)z \dd z
\end{equation}
we can also write
\begin{equation}\label{UG}
	U_G[p]=V(p)+U_0[p],
\end{equation}
thus separating the hysteretic from the non-hysteretic part.

For our purposes, we adopt the following hypothesis on the Preisach density.
\begin{hypo}\label{hpsi}
	There exists a function \(\psi^*\in L^1(0,\infty)\) such that for a.\,e. \((r,v)\in (0,\infty)\times\real\) we have \(0 \le \psi(r,v) \le \psi^*(r)\) and 
	$$
	C_\psi^* := \int_0^\infty(1+r^2)\psi^*(r) \dd r < \infty.$$
\end{hypo}
A straightforward computation shows that $G_0$ (and, consequently, $G$) are Lipschitz continuous in $C[0,T]$. Indeed, by \eqref{lip_frak} and Hypothesis \ref{hpsi} we obtain for $p_1,p_2\in C[0,T]$ and $t\in[0,T]$ that
\begin{equation}\label{lip_G}
	|G_0[p_2](t)-G_0[p_1](t)| = \bigg|\int_0^\infty\int_{\play_r[p_1](t)}^{\play_r[p_2](t)}\psi(r,v)\dd v\dd t\bigg| \le C_\psi^*\max_{\tau\in[0,t]}|p_2(\tau)-p_1(\tau)|.
\end{equation}
Moreover, the Preisach potential is continuous from $L^2(\Omega;C[0,T])$ to $L^1(\Omega;C[0,T])$, as it holds
\begin{equation}\label{cont_U0}
	|U_0[p_2](t)-U_0[p_1](t)| \le \frac{C_\psi^*}{2} \max_{\tau \in [0,t]} |p_2-p_1|(\tau) \left(|p_2(t)|+|p_1(t)|+2\right).
\end{equation}
From Hypothesis \ref{hpsi} and identity \eqref{id_frak} for the play we also obtain
\begin{equation}\label{eneg2}
	0 < U_0[p] \le C_\psi^*(1+|p|)^2, \qquad |D_0[p]_t| \le C|p_t|.
\end{equation}
The Preisach operator admits also a family of ``nonlinear'' energies. As a consequence of \eqref{id_frak}, we have for a.\,e. $t$ the inequality
$$\play_r[p]_t(p-\play_r[p])\ge 0,$$
hence
$$\play_r[p]_t(h(p)-h(\play_r[p]))\ge 0$$
for every nondecreasing function $h:\real\rightarrow\real$. Hence, for every absolutely continuous input $p$, a counterpart of \eqref{en_id0} in the form
\begin{equation}\label{Uh_ineq}
	G_0[p]_th(p)-U_h[p]_t\ge 0\text{ a.\,e.}
\end{equation}
holds with a modified potential
\begin{equation}\label{Uh}
	U_h[p](t)=\int_0^\infty\int_0^{\play_r[p](t)}h(v)\psi(r,v)\dd v\dd r.
\end{equation}
This is related to the fact that for every absolutely continuous nondecreasing function $\hat{h}:\real\rightarrow\real$, the mapping $G_{\hat{h}} := G_0\circ\hat{h}$ is also a Preisach operator, see \cite{kpr}.


\section{Statement of the problem}\label{sta}

We introduce the spaces
\be{spx}
X = W^{1,2}(\Omega)\,, \quad X_0 = \{\psi \in W^{1,2}(\Omega;\real^3): \psi\big|_{\partial\Omega} = 0\}\,, \quad X_{q^*} = W^{1,q^*}(\Omega)
\ee
for some $q^*>2$ that will be specified below in Theorem \ref{t1}. Taking into account the boundary conditions \eqref{boufull}, we consider \eqref{e1}--\eqref{e4} in variational form
\begin{align}
	&
	\begin{aligned}
		&\int_{\Omega}\big((\chi+\rho^*(1-\chi))(f(p)+G_0[p]+\dive u)\big)_t\,\phi \dd x+\int_{\Omega}\mu(p)\nabla p\cdot\nabla\phi \dd x\\
		&=\int_{\partial\Omega}\alpha(x)(p^*-p)\,\phi \dd s(x),
	\end{aligned}
	\label{e1w}
	\\[2 mm]
	&
	\begin{aligned}
		&\int_{\Omega}(P[\nabla_s u]+\BB\nabla_s u_t):\nabla_s\psi \dd x-\int_{\Omega}\big(p(\chi+\rho^*(1-\chi))+\beta(\theta-\theta_c)\big)\,\dive\psi \dd x\\
		&=\int_{\Omega}g\cdot\psi \dd x,
	\end{aligned}	\label{e2w}
	\\[2 mm]
	&
	\begin{aligned}
		&\int_{\Omega}\Bigg(\mathcal{C}_V(\theta)_t-\BB\nabla_s u_t:\nabla_s u_t-\|D_P[\nabla_s u]_t\|_*-\mu(p)|\nabla p|^2-\dch|D_0[p]_t|\\
		& \quad -\gamma(\theta,\dive u)\chi_t^2+\left(\frac{L}{\theta_c}\chi_t+\beta\dive u_t\right)\theta\Bigg)\zeta \dd x+\int_{\Omega}\kappa(\theta)\nabla\theta\cdot\nabla\zeta \dd x \\
		&=\int_{\partial\Omega}\omega(x)(\theta^*-\theta)\,\zeta \dd x,
	\end{aligned}
	\label{e3w}
	\\[2 mm]
	&\gamma(\theta,\dive u)\chi_t+\partial I_{[0,1]}(\chi)\ni(1-\rho^*)\!\left(\Phi(p)+pG_0[p]-U_0[p]+p\dive u\right)+L\left(\frac{\theta}{\theta_c}-1\right) \ \text{a.\,e.}
	\label{e4w}
\end{align}
for a.\,e. $t \in (0,T)$ and all test functions $\phi \in X$, $\psi \in X_0$ and $\zeta \in X_{q^*}$. Note that we split the capillary hysteresis terms in hysteretic and non-hysteretic part according to \eqref{G}, \eqref{U_D}--\eqref{UG}. This is done in view of the regularization performed in Section \ref{cutoff}, where only the non-hysteretic part will be affected by the cut-off.

We assume the following hypotheses to hold.

\begin{hypo}\label{h1}
There exist constants $A^\flat > 0$, $B^\flat >0$, $\bar\theta > 0$ such that
\begin{itemize}
	\item[{\rm (i)}] $\AA_e$, $\AA_h$, $\BB$ are constant symmetric positive definite fourth order tensors such that $\AA_e\,\xi:\xi \ge A^\flat |\xi|^2$, $\AA_h\,\xi:\xi \ge A^\flat |\xi|^2$, $\BB\,\xi:\xi \ge B^\flat |\xi|^2$ for all $\xi \in \real^{3\times 3}$;
	\item[{\rm (ii)}] $g \in L^\infty(0,T; L^2(\Omega;\real^3)) \cap W^{1,2}(0,T; L^2(\Omega;\real^3))$ is a given function and there exists a function $\widehat{G} \in L^4(\Omega\times(0,T))$ such that $g = -\nabla \widehat{G}$;
	\item[{\rm (iii)}] $\alpha \in W^{1,\infty}(\partial\Omega)$, $\alpha(x)\ge 0$ a.\,e. and $\int_{\partial\Omega}\alpha(x)\dd s(x)>0$; $\omega \in L^\infty(\partial\Omega)$, $\omega(x)\ge 0$ a.\,e. and $\int_{\partial\Omega}\omega(x)\dd s(x)>0$;
	\item[{\rm (iv)}] $p^* \in L^{\infty}(\partial\Omega\times(0,T))$ and $p^*_t \in L^2(\partial\Omega\times(0,T))$, $\theta^*\in L^\infty(\partial\Omega\times(0,T))$, $\theta^*_t\in L^2(\partial\Omega\times(0,T))$, $\theta^*(x,t) \ge \bar\theta$ a.\,e.;
	\item[{\rm (v)}] $p^0 \in L^\infty(\Omega) \cap W^{2,2}(\Omega)$, $u^0 \in X_0 \cap W^{1,4}(\Omega;\real^3)$, $\theta^0 \in L^\infty(\Omega) \cap W^{1,2}(\Omega)$, $\theta^0(x) \ge \bar\theta$ a.\,e., $\chi^0 \in L^{\infty}(\Omega)$, $\chi^0(x) \in [0,1]$ a.\,e.
\end{itemize}
We also assume that there exist constants $f^\sharp > f^\flat >0$, $\nu \in (0, 1/2]$, $\mu^\flat >0$, $c^\sharp > c^\flat > 0$, $1/2 \le b < \hat{b} < 1$, $\kappa^\sharp > \kappa^\flat > 0$, $0 < a < 1-b$, $a<\hat{a}<\frac{(8+3a+2b)(1+b)}{7-2b}$, $\gamma^\sharp > \gamma^\flat > 0$ such that the nonlinearities satisfy the following conditions:
\begin{itemize}
	\item[{\rm (vi)}] $G[p]=f(p)+G_0[p]$ where $f: \real \to (C_\psi^-,1-C_\psi^+)$ with $C_\psi^\pm$ from \eqref{G0b} is a continuously differentiable function, $f^\flat(1+ |p|)^{-1-\nu} \le f'(p) \le f^\sharp$ for all $p\in \real$, and $G_0$ is the Preisach operator from Section \ref{sec_prei} with density function satisfying Hypothesis \ref{hpsi} and with potential $U_0$ and dissipation $D_0$ as in \eqref{U0_D0};
	\item[{\rm (vii)}] $\mu: \real \to \real$ is a continuous function, $\mu(p) \ge \mu^\flat$ for all $p \in \real$;
	\item[{\rm (viii)}] $\mathcal{C}_V : [0,\infty) \rightarrow [0,\infty)$ is a continuously differentiable function, $\mathcal{C}_V'(\theta) =: c_V(\theta)$ is such that $c^\flat(1+\theta^{b}) \le c_V(\theta) \le c^\sharp(1+\theta^{\hat{b}})$ for all $\theta \ge 0$;
	\item[{\rm (ix)}] $\kappa : [0,\infty) \rightarrow [0,\infty)$ is a continuous function, $\kappa^\flat(1+\theta^{1+a}) \le \kappa(\theta) \le \kappa^\sharp(1+\theta^{1+\hat{a}})$	for all $\theta \ge 0$;
	\item[{\rm (x)}] $\gamma : [0,\infty) \times [0,\infty) \rightarrow [0,\infty)$ is a continuous function, $\gamma^\flat(1+\theta+|\dive u|^2) \le \gamma(\theta,\dive u) \le \gamma^\sharp(1+\theta+|\dive u|^2)$ for all $\theta\ge 0,\,u\in\real^3$;
	\item[{\rm (xi)}] $P : C([0,T];\tens) \to C([0,T];\tens)$ is the constitutive operator of elastoplasticity with dissipation operator $D_P$ defined below in \eqref{plast}--\eqref{enep}.
\end{itemize}
\end{hypo}

\begin{remark}{\em
In this remark we comment on the more technical hypotheses.
\begin{itemize}[nolistsep,noitemsep]
	\item[(vi)] The growth condition for $f$ is in agreement with the physical requirement that $f'$ has to degenerate when $p \to \pm\infty$. The specific form of the lower bound will play a substantial role in the Moser iteration argument.
	\item[(viii)] The growth condition for $c_V$ will be of fundamental importance in Subsection \ref{subsec_estu} where, in order to estimate $\dive u_t$ in $L^q(0,T;L^2(\Omega))$ with an exponent $q>4$, we will need a higher integrability (in space) for the temperature than simply $L^\infty(0,T;L^1(\Omega))$.
	\item[(ix)] The tangled bound
	$$\hat{a} < \frac{(8+3a+2b)(1+b)}{7-2b}$$
	for the growth exponent of the function $\kappa$ is required in Subsection \ref{hot}, where we apply an iterative method in order to derive higher order estimates for the temperature.
	\item[(x)] The dependence of the relaxation coefficient $\gamma$ on both $\theta$ and $\dive u$ is uncommon but crucial for obtaining estimates \eqref{esti0g} and \eqref{esti13}.
\end{itemize}
}
\end{remark}

We model the elastoplasticity following \cite{lch}.  We assume that a convex subset $0\in Z \subset \tens$ with nonempty interior representing the admissible plastic stress domain is given in the space $\tens$ of symmetric tensors, and that the constitutive relation between the strain tensor $\ve$ and the stress tensor $\sigma$ involves two fourth order tensors $\AA_h$ (the kinematic hardening tensor) and $\AA_e$ (the elasticity tensor). We define the constitutive operator $P$ by the formula
\be{plast}
P[\ve] = \AA_h \ve + \sigma^p,
\ee
where $\sigma^p$ is the solution of the variational inequality
\be{vari}
\sigma^p \in Z, \ \big(\ve_t - \AA_e^{-1}\sigma^p_t\big):\big(\sigma^p - z\big) \ge 0 \quad \mbox{a.\,e. } \ \forall z \in Z\,, \quad \sigma^p(0) = Q_Z(\ve(0))
\ee
for a given $\ve \in W^{1,1}(0,T;\tens)$, where $Q_Z: \tens \to Z$ is the orthogonal projection onto $Z$. The variational inequality \eqref{vari} has a unique solution $\sigma^p \in W^{1,1}(0,T;\tens)$ and the solution mapping
$$P : W^{1,1}(0,T;\tens) \to W^{1,1}(0,T;\tens) : \ve \mapsto \sigma^p$$
is strongly continuous, see \cite{kb}. It holds
\begin{equation}\label{ivP}
\left| P[\ve]_t \right| \le |\ve_t|.
\end{equation}

The energy potential $U_P$ and the dissipation operator $D_P$ associated with $P$ are defined by the formula
\be{enep}
U_P[\ve] = \frac12\AA_h\ve:\ve + \frac12\AA_e^{-1}\sigma^p : \sigma^p, \qquad D_P[\ve] = \ve - \AA_e^{-1}\sigma^p.
\ee
Let $M_{Z^*}$ denote the Minkowski functional of the polar set $Z^*$ to $Z$. The energy identity
\be{enep1}
P[\ve]: \ve_t - U_P[\ve]_t = \left\| D_P[\ve]_t \right\|_* \ \mbox{a.\,e.},
\ee
where $\|\cdot\|_* = M_{Z^*}(\cdot)$ is a seminorm in $\real^{3 \times 3}_{\rm sym}$, and the inequalities
\be{enep2}
U_P[\ve] \ge \frac{A^\flat}{2}|\ve|^2, \qquad \left\| D_P[\ve]_t \right\|_* \le C|\ve_t|
\ee
hold for all inputs $\ve \in W^{1,1}(0,T;\tens)$.

The operator $P$ can be extended to a continuous operator in the space $C([0,T];\tens)$ in the sense that if $\{\ve_m; m\in \nat\}$ is a sequence in $C([0,T];\tens)$, then
\be{con1}
\lim_{m\to \infty}\max_{t \in [0,T]} |\ve_m(t) - \ve(t)| = 0 \ \Longrightarrow \ \lim_{m\to \infty}\max_{t \in [0,T]}|P[\ve_m](t) - P[\ve](t)| = 0.
\ee
For two inputs $\ve_1, \ve_2 \in W^{1,1}(0,T;\tens)$ we denote $\sigma_i = P[\ve_i]$, $i=1,2$. Then
\begin{align}\label{stop1} 
(\sigma_1 - \sigma_2):\left((\ve_1)_t - (\ve_2)_t\right) &\ge \frac12\frac{\dd}{\dd t} \Big(\AA_h(\ve_1 - \ve_2):(\ve_1 - \ve_2) + \AA_e^{-1}(\sigma^p_1 - \sigma^p_2):(\sigma^p_1 - \sigma^p_2)\Big)\ \mbox{a.\,e.},\\ \label{stop2} 
|\sigma_1(t) - \sigma_2(t)| &\le C \left(|\ve_1(0) - \ve_2(0)| + \int_0^t |(\ve_1)_t - (\ve_2)_t|(\tau)\dd\tau\right) \quad \forall t\in[0,T]
\end{align}
with a constant $C$ depending only on $\AA_h$ and $\AA_e$.

For inputs $\ve \in L^2(\Omega; W^{1,1}(0,T; \tens))$ we obtain from \eqref{stop2} similarly as in \cite[Formula (6.25)]{dkr} the inequality
\be{stop3}
|\nabla \sigma(x,t)| \le C\left(|\nabla\ve(x,0)| + \int_0^t |\nabla\ve_t(x,\tau)|\dd\tau\right) \quad \mbox{a.\,e.}
\ee

The main result of the paper reads as follows.

\begin{theorem}\label{t1}
Let Hypothesis \ref{h1} hold. Then there exists a solution $(p,u,\theta,\chi)$ to the system \eqref{e1w}--\eqref{e4w} with initial conditions \eqref{ini} with the regularity
\begin{itemize}[nolistsep,noitemsep]
	\item $p \in L^\infty(\Omega\times(0,T))$, $p_t \in L^2(\Omega\times(0,T))$, $M(p) \in L^2(0,T;W^{2,2}(\Omega))$ with $M(p)$ given by \eqref{M};
	\item $u_t \in L^q(0,T;X_0 \cap W^{1,q}(\Omega;\real^3))$ for all $q<\frac{(8+3a+2b)(4+b)}{7-2b}$, $\nabla_s u \in L^2(\Omega;C([0,T];\tens))$;
	\item $\theta \in L^q(\Omega\times(0,T))$ for all $q<\frac{(8+3a+2b)(4+b)}{7-2b}$, $\nabla\theta \in L^2(\Omega\times(0,T);\real^3)$,\\ $\theta_t \in L^2(0,T;W^{-1,q^*}(\Omega))$ with $q^*>2$ given by \eqref{qstar};
	\item $\chi \in L^q(\Omega;C[0,T])$, $\chi_t \in L^q(\Omega\times(0,T))$ for all $q \in [1,\infty)$, $\chi(x,t) \in [0,1]$ a.\,e.
\end{itemize}
\end{theorem}

The reason why we do not specify the precise value of $q^*$ here is that it relies on a certain number of intermediate computations that cannot be detailed at this stage. The proof of Theorem \ref{t1} will be divided into several steps. In order to eliminate possible degeneracy of the functions $f$ and $\mu$, we start by regularizing the problem by means of a large parameter $R$. Then we prove that this regularized problem admits a solution by the standard Faedo-Galerkin method: here the parameter $R$ will be of fundamental importance in order to gain some regularity. Once we have derived suitable estimates, we pass to the limit in the Faedo-Galerkin scheme. The second part of the proof will consist in the derivation of a priori estimates independent of $R$, which will allow us to pass to the limit in the regularized system and infer the existence of a solution with the desired regularity.
\smallskip

In what follows, we denote by $C$ any positive constant depending only on the data, by $C_R$ any constant depending on the data and on $R$ and by $C_{R,\eta}$ any constant depending on the data, on $R$ and on $\eta$, all independent of the dimension $n$ of the Galerkin approximation. Furthermore, we denote by $|v|_r$ the $L^r(\Omega)$-norm of a function $v\in L^r(\Omega)$ or $v\in L^r(\Omega;\real^3)$ for $r\in [1,\infty]$, and the norm of a function $v\in W^{1,r}(\Omega)$ will be denoted by $|v|_{1;r}$. We systematically use the Korn's inequality (see \cite {nh})
\begin{equation}\label{korn}
\io |\nas w|^2(x)\dd x \ge c \|w\|_{W^{1,2}(\Omega;\real^3)}^2
\end{equation}
for every $w \in X_0$ with a constant $c>0$ independent of $w$. We will also often use the Poincar\'e inequality (see \cite{gt,lu}) in the form
\begin{equation}\label{poin}
	|v|_{1;2}^2 \le C \left(\io|\nabla v|^2(x)\dd x + \ipo \gamma(x)|v|^2(x)\dd s(x)\right)
\end{equation}
for functions $v\in W^{1,2}(\Omega)$ provided $\gamma \in L^\infty(\partial\Omega)$ is such that $\gamma \ge 0$ a.\,e. and $\ipo \gamma(x)\dd s(x) >0$. Finally, let us recall the Gagliardo-Nirenberg inequality (see \cite{bin,gt}) for $v\in W^{1,r}(\Omega)$ on a bounded Lipschitzian domain $\Omega\subset \real^N$ in the form
\begin{equation}\label{gn}
|v|_q \le C|v|_s^{1-\delta}|v|_{1;r}^\delta
\end{equation}
with $r<N$, $s < q < (rN)/(N-r)$ and with a constant $C$ depending only on $q,r,s$, where
$$\delta = \frac{\frac1s - \frac1q}{\frac1s -\frac1r + \frac1N}\,.$$


\section{Cut-off system}\label{cutoff}

We choose a regularizing parameter $R>1$, and first solve a cut-off system with the intention to let $R\to\infty$.

For $z \in \real$ we denote by
\begin{equation}\label{QR}
	Q_R(z)=\max\{-R,\min\{z,R\}\}
\end{equation}
the projection of $\real$ onto $[-R,R]$. Then we cut-off some nonlinearities by setting
\begin{equation}\label{fR}
	f_R(p) = 
	\left\{
	\begin{array}{ll}
		f(p) & \for |p|\le R\\
		f(R) + f'(R)(p-R) & \for p>R\\
		f(-R) + f'(-R)(p+R) & \for p<-R
	\end{array}
	\right.,
\end{equation}
\begin{equation}\label{R}
	\Phi_R(p)=\int_0^p f_R(z) \dd z, \qquad V_R(p)=pf_R(p)-\Phi_R(p)=\int_0^p f_R'(z)z \dd z,
\end{equation}
\begin{equation}\label{muR}
	\mu_R(p) = \mu(Q_R(p)) = \left\{
	\begin{array}{ll}
		\mu(p) & \for |p|\le R\\
		\mu(R) & \for p>R\\
		\mu(-R) & \for p<-R
	\end{array}
	\right.,
\end{equation}
\begin{equation}\label{gammaR}
	\gamma_R(p,\theta,\dive u)=\gamma(Q_R(\theta^+)+(p^2-R^2)^+,\dive u)
\end{equation}
for $p,\theta,\dive u \in \real$. Note that by Hypothesis \ref{h1}\,(vi) we deduce that $|f_R(p)| \le |f(0)|+f^\sharp|p|$, from which
\begin{equation}\label{fR_esti}
	|f_R(p)| \le C\left(1+|p|\right), \quad |\Phi_R(p)| \le C\left(1+p^2\right), \quad C\left(|p|^{1-\nu}-1\right) \le V_R(p) \le C\,p^2,
\end{equation}
and also, from Hypothesis \ref{h1}\,(x),
\begin{equation}\label{gammaR_esti}
	\begin{aligned}
		\gamma^\flat\left(1+Q_R(\theta^+)+(p^2-R^2)^++|\dive u|^2\right) &\le \gamma_R(p,\theta,\dive u) \\
		&\le \gamma^\sharp\left(1+Q_R(\theta^+)+(p^2-R^2)^++|\dive u|^2\right).
	\end{aligned}
\end{equation}
We replace \eqref{e1w}--\eqref{e4w} by the cut-off system
\begin{align}
	&
	\begin{aligned}
		&\int_{\Omega}\big((\chi+\rho^*(1-\chi))(f_R(p)+G_0[p]+\dive u)\big)_t\,\phi \dd x+\int_{\Omega}\mu_R(p)\nabla p\cdot\nabla\phi \dd x\\
		&=\int_{\partial\Omega}\alpha(x)(p^*-p)\,\phi \dd s(x),
	\end{aligned}
	\label{e1r}
	\\[2 mm]
	&
	\begin{aligned}
		&\int_{\Omega}\big((P[\nabla_s u]+\BB\nabla_s u_t):\nabla_s\psi\big) \dd x\\
		& \ -\int_{\Omega}\big(p(\chi+\rho^*(1-\chi))+\beta(Q_R(\theta^+)-\theta_c)\big)\,\dive\psi \dd x=\int_{\Omega}g\cdot\psi \dd x,
	\end{aligned}
	\label{e2r}
	\\[2 mm]
	&
	\begin{aligned}
		&\io \Bigg(\mathcal{C}_V(\theta)_t-\BB\nabla_s u_t:\nabla_s u_t-\|D_P[\nabla_s u]_t\|_*-\mu_R(p)\,Q_R(|\nabla p|^2)-\dch|D_0[p]_t| \\
		& \ -\gamma_R(p,\theta,\dive u)\chi_t^2+\left(\frac{L}{\theta_c}\chi_t+\beta\dive u_t\right)Q_R(\theta^+)\Bigg)\zeta \dd x + \int_{\Omega}\kappa(Q_R(\theta^+))\nabla\theta\cdot\nabla\zeta \dd x \\
		& = \int_{\partial\Omega}\omega(x)(\theta^*-\theta)\,\zeta \dd s(x),
	\end{aligned}
	\label{e3r}
	\\[2 mm]
	&
	\begin{aligned}
		&\gamma_R(p,\theta,\dive u)\chi_t+\partial I_{[0,1]}(\chi)\\
		&\ni(1-\rho^*)\left(\Phi_R(p)+p\,G_0[p]-U_0[p]+p\,\dive u\right)+L\left(\frac{Q_R(\theta^+)}{\theta_c}-1\right) \ \text{a.\,e.}
	\end{aligned}
\label{e4r}
\end{align}
for all test functions $\phi, \zeta \in X$ and $\psi \in X_0$. For the system \eqref{e1r}--\eqref{e4r} the following result holds true.

\begin{propo}\label{p1}
	Let Hypothesis \ref{h1} hold and let $R>1$ be given. Then there exists a solution $(p,u,\theta,\chi)$ to \eqref{e1r}--\eqref{e4r}, \eqref{ini} with the regularity
	\begin{itemize}[nolistsep,noitemsep]
		\item $p \in L^q(\Omega;C[0,T])$ for all $q \in [1,6)$, $\nabla p \in L^2(\Omega\times(0,T);\real^3)$, $p_t \in L^2(\Omega\times(0,T))$; 
		\item $u_t \in L^2(0,T;X_0)$, $\nabla_s u_t \in L^4(\Omega\times(0,T);\tens)$; 
		\item $\theta \in L^2(\Omega\times(0,T))$, $\nabla\theta \in L^\infty(0,T;L^2(\Omega;\real^3))$, $\theta_t \in L^2(\Omega\times(0,T))$;
		\item $\chi \in L^q(\Omega;C[0,T])$, $\chi_t \in L^q(\Omega\times(0,T))$ for all $q \in [1,\infty)$.
	\end{itemize}
\end{propo}

We split the proof of Proposition \ref{p1} in two steps. First, in Subsection \ref{subsec_reg}, we further regularize the system by means of a small parameter $\eta > 0$ in order to obtain some extra-regularity for the gradient of the capillary pressure. Then, in Subsection \ref{subsec_gal}, we solve this new problem by Galerkin approximations. Here the extra-regularization will be of fundamental importance in order to pass to the limit in the nonlinearity $Q_R(|\nabla p^{(n)}|^2)$, where $n$ is the dimension of the Galerkin scheme. As a last step, we let $\eta \to 0$.


\subsection{$W^{2,2}$-regularization of the capillary pressure}\label{subsec_reg}

We define the functions
\begin{equation}\label{MKR}
	M_R(p) := \int_0^p\mu_R(z)\dd z, \qquad K_R(\theta) := \int_0^\theta\kappa(Q_R(z^+))\dd z
\end{equation}
for $p,\theta \in \real$, and introduce the new variables $v=M_R(p)$, $z=K_R(\theta)$. We then choose another regularizing parameter $\eta \in (0,1)$ and consider the following system in the unknowns $v,u,z,\chi$:
\begin{align}
	&
	\begin{aligned}
		&\io \big((\chi+\rho^*(1-\chi))(f_R(M_R^{-1}(v))+G_0[M_R^{-1}(v)]+\dive u)\big)_t\,\phi \dd x\\
		& + \io \left(\nabla v\cdot\nabla\phi+\eta\Delta v\Delta\phi\right) \dd x = \int_{\partial\Omega} \alpha(x)(p^*-M_R^{-1}(v))\,\phi \dd s(x),
	\end{aligned}
	\label{e1e}
	\\[2 mm]
	&
	\begin{aligned}
		&\io \big((P[\nabla_s u]+\BB\nabla_s u_t):\nabla_s\psi\big) \dd x\\
		& \ - \io \big(M_R^{-1}(v)(\chi+\rho^*(1-\chi))+\beta(Q_R((K_R^{-1}(z))^+)-\theta_c)\big)\,\dive\psi \dd x = \io g\cdot\psi \dd x,
	\end{aligned}
	\label{e2e}
	\\[2 mm]
	&
	\begin{aligned}
		&\io \Bigg(\mathcal{C}_V(K_R^{-1}(z))_t-\BB\nabla_s u_t:\nabla_s u_t-\|D_P[\nabla_s u]_t\|_*-\mu_R(M_R^{-1}(v))Q_R(|\nabla(M_R^{-1}(v))|^2) \\
		& \ -(\chi+\rho^*(1-\chi))|D_0[M_R^{-1}(v)]_t|-\gamma_R(M_R^{-1}(v),K_R^{-1}(z),\dive u)\chi_t^2 \\
		& + \left(\frac{L}{\theta_c}\chi_t+\beta\dive u_t\right)Q_R((K_R^{-1}(z))^+)\Bigg)\zeta \dd x + \int_{\Omega}\nabla z\cdot\nabla\zeta \dd x \\
		&= \int_{\partial\Omega} \omega(x)(\theta^*-K_R^{-1}(z))\,\zeta \dd s(x),
	\end{aligned}
	\label{e3e}
	\\[2 mm]
	&
	\begin{aligned}
		&\gamma_R(M_R^{-1}(v),K_R^{-1}(z),\dive u)\chi_t+\partial I_{[0,1]}(\chi)\\
		&\ni(1-\rho^*)\left(\Phi_R(M_R^{-1}(v))+M_R^{-1}(v)\,G_0[M_R^{-1}(v)]-U_0[M_R^{-1}(v)]+M_R^{-1}(v)\,\dive u\right) \\
		& + L\left(\frac{Q_R((K_R^{-1}(z))^+)}{\theta_c}-1\right) \ \text{a.\,e.}
	\end{aligned}
	\label{e4e}
\end{align}
with test functions $\phi \in W^{2,2}(\Omega)$, $\psi \in X_0$ and $\zeta \in X$.


\subsection{Galerkin approximations}\label{subsec_gal}

For each fixed $R>1$ and $\eta\in (0,1)$, system \eqref{e1e}--\eqref{e4e} will be solved by Faedo-Galerkin approximations. To this end, let $\mathcal{W}=\left\{e_i : i=0,1,2,\dots\right\}\subset L^2(\Omega)$ be the complete orthonormal systems of eigenfunctions defined by
$$
	-\Delta e_i=\lambda_i e_i\quad \text{in }\Omega,\qquad\nabla e_i\cdot n\big|_{\partial\Omega}=0
$$
with $\lambda_0=0$, $\lambda_i>0$ for $i\ge 1$. Given $n\in\nat$, we approximate $v$ and $z$ by the finite sums
$$v^{(n)}(x,t)=\sum_{i=0}^{n}v_i(t)e_i(x), \quad z^{(n)}(x,t)=\sum_{k=0}^{n}z_k(t)e_k(x)$$
where the coefficients $v_i,z_k:[0,T]\to\real$ and $u^{(n)},\chi^{(n)}$ will be determined as the solution of the system
\begin{align}
	&
	\begin{aligned}
		&\int_\Omega\big((\chi^{(n)}+\rho^*(1-\chi^{(n)}))(f_R(p^{(n)})+G_0[p^{(n)}]+\dive u^{(n)})\big)_t\,e_i\dd x+\left(\lambda_i+\eta\lambda_i^2\right) v_i\\
		&=\int_{\partial\Omega}\alpha(x)(p^*-p^{(n)})\,e_i\dd s(x),
	\end{aligned}
	\label{e1g}
	\\[2 mm]
	&
	\begin{aligned}
		&\int_\Omega \big((P[\nabla_s u^{(n)}]+\BB\nabla_s u^{(n)}_t):\nabla_s\psi \dd x\\
		&-\int_\Omega\Big(p^{(n)}(\chi^{(n)}+\rho^*(1-\chi^{(n)}))+\beta(Q_R((\theta^{(n)})^+)-\theta_c)\Big)\,\dive\psi \dd x=\int_\Omega g\cdot\psi \dd x,
	\end{aligned}
	\label{e2g}
	\\[2 mm]
	&
	\begin{aligned}
		&\int_\Omega\Bigg(\mathcal{C}_V(\theta^{(n)})_t-\BB\nabla_s u^{(n)}_t:\nabla_s u^{(n)}_t-\|D_P[\nabla_s u^{(n)}]_t\|_*-\mu_R(p^{(n)})Q_R(|\nabla p^{(n)}|^2) \\
		& \ -(\chi^{(n)}+\rho^*(1-\chi^{(n)}))|D_0[p^{(n)}]_t|-\gamma_R(p^{(n)},\theta^{(n)},\dive u^{(n)})|\chi^{(n)}_t|^2 \\
		& \ +\left(\frac{L}{\theta_c}\chi^{(n)}_t+\beta\dive u^{(n)}_t\right)Q_R((\theta^{(n)})^+)\Bigg)\,e_k\dd x + \lambda_k z_k=\int_{\partial\Omega}\omega(x)(\theta^*-\theta^{(n)})\,e_k\dd s(x),
	\end{aligned}
	\label{e3g}
	\\[2 mm]
	&
	\begin{aligned}
		&\gamma_R(p^{(n)},\theta^{(n)},\dive u^{(n)})\chi^{(n)}_t+\partial I_{[0,1]}(\chi^{(n)})\\
		&\ni(1-\rho^*)\!\left(\Phi_R(p^{(n)})+p^{(n)}G_0[p^{(n)}]-U_0[p^{(n)}]+p^{(n)}\,\dive u^{(n)}\right)+L\left(\frac{Q_R((\theta^{(n)})^+)}{\theta_c}-1\right) \ \text{a.\,e.}
	\end{aligned}
	\label{e4g}
\end{align}
for $i,k=0,1,\dots,n$ and for all $\psi \in X_0$, and with $p^{(n)}:=M_R^{-1}(v^{(n)})$, $\theta^{(n)}:=K_R^{-1}(z^{(n)})$. We prescribe the initial conditions
\begin{equation}\label{inig}
	\left.
	\begin{array}{rcl}
		v_i(0) &=& \int_\Omega M_R(p^0(x))\,e_i(x)\dd x,\\ [2 mm]
		u^{(n)}(x,0) &=& u^0(x),\\ [2 mm]
		z_k(0) &=& \int_\Omega K_R(\theta^0(x))\,e_k(x)\dd x,\\ [2 mm]
		\chi^{(n)}(x,0) &=& \chi^0(x).
	\end{array}
	\right\rbrace
\end{equation}

This is an ODE system coupled with a nonlinear PDE \eqref{e2g}. It is nontrivial to prove that such a system admits a unique strong solution. We proceed as follows. For a given function $w \in L^r(\Omega\times(0,T))$ consider the equation
\begin{equation}\label{ex1}
	\io \BB\nabla_s u_t(x,t):\nabla_s\psi(x) \dd x + \io P[\nabla_s u](x,t):\nabla_s\psi(x) \dd x = \io w(x,t)\,\dive\psi(x) \dd x,
\end{equation}
which is to be satisfied for every $\psi \in X_0$ a.\,e. in $(0,T)$ together with an initial condition $u(x,0)=u^0(x)$, $u^0 \in X_0 \cap W^{1,r}(\Omega;\real^n)$ and boundary condition $u=0$ on $\partial\Omega$.

\textbf{Step 1.} As an application of the $L^r$-regularity theory for elliptic systems in divergence form (see e.\,g. \cite{sima}) we can conclude that for every $w \in L^r(\Omega\times(0,T))$ with some $r \in [2,\infty)$ the problem
$$\io \BB\nabla_s u_t(x,t):\nabla_s\psi(x) \dd x = \io w(x,t)\,\dive\psi(x) \dd x$$
has for a.\,e. $t\in (0,T)$ a unique solution such that $\nabla_s u_t(\cdot, t) \in L^r(\Omega;\tens)$, and it holds
$$
\io |\nabla_s u_t|^r(x,t) \dd x \le C \io |w|^r(x,t) \dd x \quad \text{for a.\,e. } t \in (0,T).
$$
Integrating over $t \in (0,T)$ we obtain that $\nabla_s u_t \in L^r(\Omega\times (0,T);\tens)$.

\textbf{Step 2.} Let us define the convex and closed subset $\Ucal_r := \{w \in L^r(0,T;X_0) : \nabla_s w_t \in L^r(\Omega\times(0,T);\tens), \ w(x,0)=u^0(x) \text{ a.\,e.}\} \subset L^r(0,T;X_0\cap W^{1,r}(\Omega;\real^3))$. Note that the trace operator is well defined on this space, so that the initial condition makes sense. Let $\hat{u} \in \Ucal_r$, and let $u$ be the solution of the equation
$$\io P[\nabla_s \hat{u}](x,t):\nabla_s\psi(x) \dd x + \io \BB\nabla_s u_t(x,t):\nabla_s\psi(x) \dd x = \io w(x,t)\,\dive\psi(x) \dd x,$$
the existence of which follows from Step 1. We prove that the mapping $\hat{u}_t \mapsto u_t$ is a contraction with respect to a suitable norm.

Indeed, let $\hat{u}_1,\hat{u}_2$ be given, and let $u_1,u_2$ be the corresponding solutions. The difference $\bar{u}=u_1-u_2$ is the solution of the equation
$$\io \BB\nabla_s\bar{u}_t(x,t):\nabla_s\psi(x) \dd x = -\io (P[\nabla_s\hat{u}_1]-P[\nabla_s\hat{u}_2])(x,t):\nabla_s\psi(x) \dd x.$$
According to Step 1, we have
\begin{equation}\label{ex2}
	\io |\nabla_s\bar{u}_t|^r(x,t) \dd x \le C \io |P[\nabla_s\hat{u}_1]-P[\nabla_s\hat{u}_2]|^r(x,t) \dd x \quad \text{a.\,e.}
\end{equation}
By inequality \eqref{stop2} we have for a.\,e. $(x,t)$
$$|P[\nabla_s\hat{u}_1]-P[\nabla_s\hat{u}_2]|(x,t) \le C \int_0^t |\nabla_s(\hat{u}_1-\hat{u}_2)_t| (x,\tau) \dd\tau$$
with a constant $C>0$. Hence, by H\"older's inequality,
\begin{align}
	\io |P[\nabla_s\hat{u}_1]-P[\nabla_s\hat{u}_2]|^r(x,t) \dd x &\le C \io \left(\int_0^t |\nabla_s(\hat{u}_1-\hat{u}_2)_t| (x,\tau) \dd\tau\right)^r \dd x \cr
	&\le Ct^{r-1} \int_0^t\io |\nabla_s(\hat{u}_1-\hat{u}_2)_t|^r (x,\tau) \dd x\dd\tau. \label{ex3}
\end{align}
Now, set
$$W(t) = \io |\nabla_s(u_1-u_2)_t|^r (x,t) \dd x, \qquad \hat{W}(t) = \io |\nabla_s(\hat{u}_1-\hat{u}_2)_t|^r (x,t) \dd x.$$
It follows from \eqref{ex2} and \eqref{ex3} that
$$W(t) \le Ct^{r-1} \int_0^t \hat{W}(\tau) \dd\tau.$$
We now multiply both sides of the above inequality by $\e^{-Ct^r}$, and after an integration over $t \in [0,T]$ we obtain from the Fubini Theorem
$$\int_0^T \e^{-Ct^r} W(t) \dd t \le \frac1p \int_0^T \left(\e^{-C\tau^r}-\e^{-CT^r}\right)\hat{W}(\tau) \dd\tau \le \frac1p \int_0^T \e^{-Ct^r}\hat{W}(t) \dd t.$$
This means that the mapping $\hat{u}_t \mapsto u_t$ is a contraction in $L^r(0,T;X_0\cap W^{1,r}(\Omega;\real^3))$ with respect to the weighted norm
$$||| u_t ||| = \left(\int_0^T \e^{-Ct^r} \io |\nabla_s u_t|^r(x,t)\dd x\dd t\right)^{1/r},$$
hence it has a unique fixed point which is a solution of \eqref{ex1}.

\textbf{Step 3.} The mapping which with a right-hand side $w \in L^r(\Omega\times(0,T))$ associates the solution $u_t \in L^r(0,T;X_0\cap W^{1,r}(\Omega;\real^3))$ of \eqref{ex1} is Lipschitz continuous. Indeed, consider $w_1,w_2$ and the corresponding solutions $u_1,u_2$, and set as before $\bar{w}=w_1-w_2$, $\bar{u}=u_1-u_2$. As a counterpart of \eqref{ex2} we get
$$\io |\nabla_s\bar{u}_t|^r(x,t) \dd x \le C \io \left(|P[\nabla_s u_1]-P[\nabla_s u_2]|^r+|\bar{w}|^r\right)(x,t) \dd x \quad \text{a.\,e.},$$
and the computations as in \eqref{ex3} yield
$$\io |\nabla_s\bar{u}_t|^r(x,t) \dd x \le Ct^{r-1} \int_0^t\io |\nabla_s\bar{u}_t|^r(x,\tau) \dd x \dd\tau + C \io |\bar{w}|^r(x,t) \dd x \quad \text{a.\,e.}$$
We obtain the Lipschitz continuity result when we test by $\e^{-\frac{C}{p}t^r}$ and integrate over $t \in [0,T]$, similarly as in Step 2.

\smallskip

Now, coming back to our equation \eqref{e2g}, we see that it is of the form \eqref{ex1} with $w(x,t) = w^{(n)}(x,t) := p^{(n)}(\chi^{(n)}+\rho^*(1-\chi^{(n)}))+\beta(Q_R((\theta^{(n)})^+)-\theta_c)(x,t) + \widehat{G}(x,t)$, where $\widehat{G}$ is from Hypothesis \ref{h1}\,(ii). Therefore, denoting by $\mathcal{S} : w^{(n)} \mapsto u^{(n)}_t$ its associated solution operator, we conclude that \eqref{e1g}--\eqref{e4g} give rise to a system of ODEs with a locally Lipschitz continuous right-hand side containing the operator $\mathcal{S}$. 

Thus system \eqref{e1g}--\eqref{e4g} has a unique strong solution in a maximal interval of existence $[0,T_n] \subset [0,T]$. This interval coincides with the whole $[0,T]$, provided we prove that the solution remains bounded in $[0,T_n)$.

We now derive a series of estimates. Note that we decompose the auxiliary variables $v$ and $z$ instead of $p$ and $\theta$ into a Fourier series with respect to the basis $\Wcal$ because we are going to test equations \eqref{e1g} and \eqref{e3g} by nonlinear expressions of $p$ and $\theta$, namely, by their Kirchhoff transforms \eqref{MKR}. Indeed, the Galerkin method allows only to test by linear functions and their derivatives.

Moreover, we do not discretize the momentum balance equation because considering the full PDE is the only way to deduce compactness of the sequence $\{\nabla_s u^{(n)}_t\}$, which is needed in order to pass to the limit in some nonlinear terms. Indeed, we will not be able to control higher derivatives of $u^{(n)}$, and this will prevent us from applying the usual embedding theorems.


\subsection{Estimates independent of $n$}\label{subsec_estn}

\textbf{Estimate 1.} We test \eqref{e1g} by $v_i$ and sum up over $i=0,1,\dots,n$, and \eqref{e2g} by $\psi=u^{(n)}_t$. Then we sum up the two equations to obtain
\begin{equation}\label{eg1}
	\begin{aligned}
		&\int_\Omega(\chi^{(n)}+\rho^*(1-\chi^{(n)}))f_R(p^{(n)})_tM_R(p^{(n)})\dd x\\
		&+\int_\Omega(\chi^{(n)}+\rho^*(1-\chi^{(n)}))G_0[p^{(n)}]_tM_R(p^{(n)})\dd x\\
		&+\int_\Omega \left(\BB\nabla_s u^{(n)}_t:\nabla_s u^{(n)}_t+U_P[\nabla_s u^{(n)}]_t+\|D_P[\nabla_s u^{(n)}]_t\|_*\right) \dd x\\
		&+\int_\Omega\left(|\nabla v^{(n)}|^2+\eta\,|\Delta v^{(n)}|^2\right)\dd x+\int_{\partial\Omega}\alpha(x)(p^{(n)}-p^*)M_R(p^{(n)})\dd s(x)\\
		&=-\int_\Omega(1-\rho^*)\chi^{(n)}_t\left(f_R(p^{(n)})+G_0[p^{(n)}]+\dive u^{(n)}\right)M_R(p^{(n)})\dd x\\
		&+\int_\Omega (\chi^{(n)}+\rho^*(1-\chi^{(n)}))\,\dive u^{(n)}_t \left(p^{(n)}-M_R(p^{(n)})\right)\dd x\\
		&+\int_\Omega\beta(Q_R((\theta^{(n)})^+)-\theta_c)\,\dive u^{(n)}_t\dd x+\int_\Omega g\cdot u^{(n)}_t\dd x.
	\end{aligned}
\end{equation}
where we exploited also the energy identity \eqref{enep1}. We now define
$$V_{M,R}(p) := \int_0^p f_R'(z)M_R(z)\dd z$$
so that
\begin{align*}
	&\int_\Omega(\chi^{(n)}+\rho^*(1-\chi^{(n)}))f_R(p^{(n)})_tM_R(p^{(n)})\dd x \\
	& = \frac{\dd}{\dd t} \io (\chi^{(n)}+\rho^*(1-\chi^{(n)}))V_{M,R}(p^{(n)}) \dd x - \io (1-\rho^*)\chi^{(n)}_tV_{M,R}(p^{(n)})\dd x,
\end{align*}
and introduce the modified Preisach potential as a counterpart to \eqref{Uh}
$$U_{M,R}[p] := \int_0^\infty\int_0^{\play_r[p]}M_R(v)\,\psi(r,v)\dd v\dd r > 0$$
which satisfies
$$G_0[p]_t\,M_R(p)-U_{M,R}[p]_t \ge 0 \ \text{ a.\,e.}$$
according to \eqref{Uh_ineq}. Note that \eqref{fR} and \eqref{muR} together with Hypothesis \ref{h1}\,(vi) and (vii) yield
\begin{equation}\label{FMR}
	c_R\,p^2 \le V_{M,R}(p) \le C_R\,p^2
\end{equation}
for all $p \in \real$, with some positive constants $c_R,C_R$ depending only on $R$. Moreover, the estimate
\be{UMR}
	U_{M,R}[p] \le C_R\left(1+|p|\right)^2
\ee
holds as a counterpart of \eqref{eneg2}. 
By the definition of $v^{(n)}$ and Hypothesis \ref{h1}\,(vii) we deduce
\begin{equation}\label{nablav}
	\int_\Omega|\nabla v^{(n)}|^2\dd x = \int_\Omega|\mu_R(p^{(n)})|^2|\nabla p^{(n)}|^2\dd x \ge (\mu^\flat)^2\int_\Omega|\nabla p^{(n)}|^2\dd x.
\end{equation}
Moreover, thanks again to Hypothesis \ref{h1}\,(vii), the boundary term is such that
\begin{align*}
	&\int_{\partial\Omega}\alpha(x)(p^{(n)}-p^*)M_R(p^{(n)})\dd s(x) \\
	&=\int_{\partial\Omega}\alpha(x)|p^{(n)}|^2 M_R^*(p^{(n)})\dd s(x) - \int_{\partial\Omega}\alpha(x)p^*p^{(n)}M_R^*(p^{(n)})\dd s(x) \\
	&\ge \mu^\flat\int_{\partial\Omega}\alpha(x)|p^{(n)}|^2\dd s(x)- C_R\int_{\partial\Omega}\alpha(x)|p^*p^{(n)}|\dd s(x),
\end{align*}
where for $p \in \real$ we set
$$
M_R^*(p) :=
\begin{cases}
	M_R(p)/p &\text{ for }p \ne 0,\\
	M_R'(0) &\text{ for }p = 0.
\end{cases}
$$
Young's inequality and Hypothesis \ref{h1}\,(iii),\,(iv) give
\begin{equation}\label{bou1}
	\int_{\partial\Omega}\alpha(x)(p^{(n)}-p^*)M_R(p^{(n)})\dd s(x) \ge \frac{\mu^\flat}{2}\int_{\partial\Omega}\alpha(x)|p^{(n)}|^2\dd s(x) - C_R.
\end{equation}
Moreover, by H\"older's inequality and Hypothesis \ref{h1}\,(ii),
\begin{align}\label{eg3}
	\io (g\cdot u^{(n)}_t)(x,t) \dd x &\le C\left(\io |u^{(n)}_t|^2(x,t) \dd x\right)^{1/2} \nonumber\\
	&\le \frac{2C}{\sqrt{c\,B^\flat}}\,\frac{\sqrt{B^\flat}}{2}\left(\io |\nabla_s u^{(n)}_t|^2 \dd x\right)^{1/2} \le C + \frac{B^\flat}{8} \io |\nabla_s u^{(n)}_t|^2 \dd x
\end{align}
where in the last line we used first Korn's inequality \eqref{korn} and then Young's inequality. Neglecting some lower order positive terms on the left-hand side, exploiting estimates \eqref{fR_esti}, \eqref{FMR} and \eqref{UMR}, and the fact that $\rho^* \le (\chi+\rho^*(1-\chi)) \le 1$ for all $\chi \in [0,1]$, from \eqref{eg1} and the subsequent computations we obtain
\begin{equation}\label{eg4}
	\begin{aligned}
		&\frac{\dd}{\dd t} \io\Big((\chi^{(n)}+\rho^*(1-\chi^{(n)}))(V_{M,R}(p^{(n)})+U_{M,R}[p^{(n)}])+U_P[\nabla_s u^{(n)}]\Big) \dd x \\
		& + \int_\Omega\bigg(|\nabla p^{(n)}|^2+\eta\,|\Delta v^{(n)}|^2+\frac{B^\flat}{2}|\nabla_s u^{(n)}_t|^2\bigg)\dd x + \int_{\partial\Omega}\alpha(x)|p^{(n)}|^2 \dd s(x) \\
		&\le C_R \left(1 + \io\left(|\chi^{(n)}_t||p^{(n)}|^2+|\chi^{(n)}_t||\dive u^{(n)}||p^{(n)}|+|p^{(n)}|^2\right)\dd x\right),
	\end{aligned}
\end{equation}
where we used also Hypothesis \ref{h1}\,(i), Young's inequality and the pointwise inequality
\begin{equation}\label{div}
	|\dive u|^2(x,t) \le 3 \, |\nabla_s u|^2(x,t) \quad\text{for all }(x,t) \in \Omega\times(0,T)
\end{equation}
to absorb $|\dive u^{(n)}_t|$ on the left-hand side together with the term coming from \eqref{eg3}.

We need to control $|\chi^{(n)}_t|$. To this aim note that \eqref{e4g} is of standard form, namely,
$$\chi^{(n)}_t + \partial I_{[0,1]}(\chi^{(n)}) \ni F^{(n)}$$
with
\begin{equation}\label{Fn}
	F^{(n)} = \frac{(1-\rho^*)\left(\Phi_R(p^{(n)})+p^{(n)} G_0[p^{(n)}]- U_0[p^{(n)}]+p^{(n)}\dive u^{(n)}\right)+L\left(Q_R((\theta^{(n)})^+)/\theta_c-1\right)}{\gamma_R(p^{(n)},\theta^{(n)},\dive u^{(n)})}\,,
\end{equation}
or, equivalently,
$$\chi^{(n)} \in [0,1], \quad (F^{(n)}-\chi^{(n)}_t)(\chi^{(n)}-\tilde{\chi}) \ge 0 \quad \text{a.\,e.} \ \ \forall \tilde{\chi} \in [0,1].$$
This yields, thanks to \eqref{eneg2}, \eqref{fR_esti} and \eqref{gammaR_esti},
\begin{equation}\label{esti0g}
	|\chi^{(n)}_t(x,t)| \le \frac{C\,(1+|p^{(n)}|^2)+|\dive u^{(n)}|^2/2+L\left|Q_R((\theta^{(n)})^+)/\theta_c-1\right|}{\gamma^\flat\left(1+(Q_R((\theta^{(n)})^+)+(|p^{(n)}|^2-R^2)^++|\dive u^{(n)}|^2\right)} \le C_R
\end{equation}
for a.\,e. $(x,t) \in \Omega\times(0,T_n)$. We now come back to \eqref{eg4} and integrate in time $\int_0^\tau \dd t$ for some $\tau \in [0,T_n]$. The initial conditions are kept under control thanks to \eqref{vari}, \eqref{enep}, \eqref{FMR}, \eqref{UMR} and Hypothesis \ref{h1}\,(v). Hence Young's inequality yields
\begin{align*}
	&\io \left(|p^{(n)}|^2+|\nabla_s u^{(n)}|^2\right)(x,\tau) \dd x + \int_0^\tau\int_\Omega\left(|\nabla p^{(n)}|^2+\eta\,|\Delta v^{(n)}|^2+|\nabla_s u^{(n)}_t|^2\right)(x,t)\dd x\dd t \\
	& + \int_0^\tau\int_{\partial\Omega}\alpha(x)|p^{(n)}|^2(x,t)\dd s(x)\dd t \le C_R \left(1 + \int_0^\tau\io \left(|p^{(n)}|^2+|\dive u^{(n)}|^2\right)(x,t) \dd x\dd t\right).
\end{align*}
Using \eqref{div} and Gr\"onwall's lemma, we thus obtain
\begin{align}
	\supess_{\tau \in (0,T_n)} \io \left(|p^{(n)}|^2 + |\nabla_s u^{(n)}|^2\right)(x,\tau)\dd x &\le C_R, \label{esti1g} \\
	\int_0^{T_n}\!\!\left(\io \Big(|\nabla p^{(n)}|^2 + |\nabla_s u^{(n)}_t|^2\Big)(x,t) \dd x + \int_{\partial\Omega} \alpha(x)|p^{(n)}|^2(x,t) \dd s(x)\right)\dd t &\le C_R, \label{esti2g}
\end{align}
and also
\begin{equation}\label{laplvn}
	\int_0^{T_n}\io |\Delta v^{(n)}|^2(x,t) \dd x \dd t \le \frac{C_R}{\eta}\,.
\end{equation}
Now, in terms of the variable $v^{(n)}=M_R(p^{(n)})$, the boundary condition is nonlinear. By the spatial $W^{2,2}$-regularity result for parabolic equations with nonlinear boundary conditions on $C^{1,1}$ domains stated and proved in \cite[Theorem 4.1]{kp}, we see that
\begin{equation}\label{esti8g} 
	\|M_R(p^{(n)})\|^2_{L^2(0,T_n;W^{2,2}(\Omega))} \le \frac{C_R}{\eta}\,.
\end{equation}
This, together with the Sobolev embedding $W^{1,2}(\Omega) \hookrightarrow L^6(\Omega)$, yields
\begin{equation*}
	\int_0^T \left(\io |\nabla p^{(n)}|^6(x,t) \dd x\right)^{1/3} \dd t \le \frac{C_R}{\eta}\,,
\end{equation*}
and since $W^{1,6}(\Omega) \hookrightarrow C(\bar{\Omega})$ we also get
\begin{equation}\label{esti10g}
	\int_0^T \supess_{x \in \Omega} |p^{(n)}|^2(x,t) \dd t \le \frac{C_R}{\eta}\,.
\end{equation}

\textbf{Estimate 2.} We test \eqref{e1g} by $\dot v_i$ and sum up over $i=0,1,\dots,n$. We get
\begin{equation}\label{eg5}
	\begin{aligned}
		&\int_\Omega\big((\chi^{(n)}+\rho^*(1-\chi^{(n)}))(f_R(p^{(n)})+G_0[p^{(n)}]+\dive u^{(n)})\big)_t\,M_R(p^{(n)})_t\dd x\\
		&+\int_\Omega\left(\nabla v^{(n)}\cdot\nabla v^{(n)}_t+\eta\,\Delta v^{(n)}\Delta v^{(n)}_t\right) \dd x=\int_{\partial\Omega}\alpha(x)(p^*-p^{(n)})M_R(p^{(n)})_t\dd s(x).
	\end{aligned}
\end{equation}
Defining
$$\hat\mu_R(p) := \int_0^p M_R'(z)z\dd z = \int_0^p \mu_R(z)z\dd z$$
for $p \in \real$, we can rewrite
\begin{align*}
	&\int_{\partial\Omega}\alpha(x)(p^{(n)}-p^*)M_R(p^{(n)})_t\dd s(x)
	\\
	&=\int_{\partial\Omega}\alpha(x)\big(\hat\mu_R(p^{(n)})-p^*M_R(p^{(n)})\big)_t\dd s(x)+\int_{\partial\Omega}\alpha(x)p^*_t M_R(p^{(n)})\dd s(x).
\end{align*}
Hence, computing the time derivative in the first summand and rearranging the terms, we can rewrite \eqref{eg5} as
\begin{equation}\label{eg6}
	\begin{aligned}
		&\frac{\dd}{\dd t}\left(\int_\Omega\frac{1}{2}\left(|\nabla v^{(n)}|^2+\eta\,|\Delta v^{(n)}|^2\right)\dd x+\int_{\partial\Omega}\alpha(x)(\hat\mu_R(p^{(n)})-p^*M_R(p^{(n)}))\dd s(x)\right)\\
		&+\int_\Omega(\chi^{(n)}+\rho^*(1-\chi^{(n)}))\left(f_R(p^{(n)})_t+G[p^{(n)}]_t\right)M_R(p^{(n)})_t\dd x\\
		&=-\int_\Omega(1-\rho^*)\chi^{(n)}_t\left(f_R(p^{(n)})+G_0[p^{(n)}]+\dive u^{(n)}\right)M_R(p^{(n)})_t\dd x\\
		&-\int_\Omega(\chi^{(n)}+\rho^*(1-\chi^{(n)}))\,\dive u^{(n)}_tM_R(p^{(n)})_t\dd x-\int_{\partial\Omega}\alpha(x)p^*_tM_R(p^{(n)})\dd s(x).
	\end{aligned}
\end{equation}
Combining \eqref{G0} with the identity \eqref{id_frak} for the play, we see that it holds
$$G_0[p^{(n)}]_t\,M_R(p^{(n)})_t = G_0[p^{(n)}]_t\,p^{(n)}_t\mu_R(p^{(n)}) \ge 0.$$
Hence by \eqref{fR}, \eqref{muR} and Hypothesis \ref{h1}\,(vi),\,(vii) we obtain the pointwise lower bound
$$(\chi^{(n)}+\rho^*(1-\chi^{(n)}))\left(f_R(p^{(n)})_t+G_0[p^{(n)}]_t\right)M_R(p^{(n)})_t \ge \rho^*f_R'(p^{(n)})\,\mu_R(p^{(n)})\,|p^{(n)}_t|^2 \ge C_R\,|p^{(n)}_t|^2.$$
We now integrate \eqref{eg6} in time $\int_0^\tau \dd t$ for some $\tau \in [0,T_n)$. Note that $\hat\mu_R(p) \ge \mu^\flat\,p^2/2$ for all $p \in \real$. Hence, arguing as for estimate \eqref{bou1} with $p^{(n)}M_R(p^{(n)})$ replaced by $\hat\mu_R(p^{(n)})$, we obtain that the boundary term on the left-hand side is such that
\begin{equation}\label{eg7}
	\int_{\partial\Omega}\alpha(x)(\hat\mu_R(p^{(n)})-p^*M_R(p^{(n)}))(x,\tau)\dd s(x)\ge \frac{\mu^\flat}{4}\int_{\partial\Omega}\alpha(x)|p^{(n)}|^2(x,\tau)\dd s(x)-C_R.
\end{equation}
Concerning the initial conditions, we employ Hypothesis \ref{h1}\,(iv) and (v). Thus, exploiting also \eqref{nablav} and \eqref{esti0g}, we get
\begin{align*}
	&\io \left(|\nabla p^{(n)}|^2+\eta\,|\Delta v^{(n)}|^2 \right)(x,\tau)\dd x + \int_{\partial\Omega}\alpha(x)|p^{(n)}|^2(x,\tau)\dd s(x) + \int_0^\tau\io |p^{(n)}_t|^2(x,t) \dd x\dd t \\
	&\le C_R\,\bigg(1 + \int_0^\tau\io \left(|\dive u^{(n)}||p^{(n)}_t|+|p^{(n)}||p^{(n)}_t|+|\dive u^{(n)}_t||p^{(n)}_t|\right) \dd x\dd t \\
	& \quad + \int_0^\tau\int_{\partial\Omega}\alpha(x)|p^*_t||p^{(n)}|\dd s(x) \dd t\bigg).
\end{align*}
Young's inequality, Hypothesis \ref{h1}\,(iii) and (iv), and estimates \eqref{div}, \eqref{esti1g}, \eqref{esti2g} give
\begin{align}
	\supess_{\tau \in (0,T_n)} \left(\io |\nabla p^{(n)}|^2(x,\tau) \dd x + \int_{\partial\Omega} \alpha(x)|p^{(n)}|^2(x,\tau)\dd s(x)\right) &\le C_R, \label{esti3g} \\
	\supess_{\tau \in (0,T_n)} \io |\Delta v^{(n)}|^2(x,\tau) \dd x &\le \frac{C_R}{\eta}, \label{esti9g} \\
	\int_0^{T_n} \io |p^{(n)}_t|^2(x,t) \dd x \dd t &\le C_R. \label{esti4g}
\end{align}

\textbf{Estimate 3.} We test \eqref{e3g} by $\dot z_k$ and sum over $k=0,1,\dots,n$. We obtain
\begin{equation}\label{eg8}
	\begin{aligned}
		&\int_\Omega\Bigg(\mathcal{C}_V(\theta^{(n)})_t-\BB\nabla_s u^{(n)}_t:\nabla_s u^{(n)}_t-\|D_P[\nabla_s u^{(n)}]_t\|_*-\mu_R(p^{(n)})Q_R(|\nabla p^{(n)}|^2) \\
		& \quad -(\chi^{(n)}+\rho^*(1-\chi^{(n)}))|D_0[p^{(n)}]_t|-\gamma_R(p^{(n)},\theta^{(n)},\dive u^{(n)})|\chi^{(n)}_t|^2 \\
		& \quad +\bigg(\frac{L}{\theta_c}\chi^{(n)}_t+\beta\dive u^{(n)}_t\bigg)Q_R((\theta^{(n)})^+)\Bigg)K_R(\theta^{(n)})_t\dd x+\int_\Omega\nabla z^{(n)}\cdot\nabla z^{(n)}_t\dd x \\
		&=\int_{\partial\Omega}\omega(x)(\theta^*-\theta^{(n)})K_R(\theta^{(n)})_t\dd s(x).
	\end{aligned}
\end{equation}
Defining
$$\hat\kappa_R(\theta) := \int_0^\theta K_R'(z)z\dd z = \int_0^\theta \kappa(Q_R(z^+))z\dd z$$
for $\theta \in \real$, we can rewrite
\begin{align*}
	&\int_{\partial\Omega}\omega(x)(\theta^{(n)}-\theta^*)K_R(\theta^{(n)})_t\dd s(x)
	\\
	&=\int_{\partial\Omega}\omega(x)\big(\hat\kappa_R(\theta^{(n)})-\theta^*K_R(\theta^{(n)})\big)_t\dd s(x)+\int_{\partial\Omega}\omega(x)\theta^*_t K_R(\theta^{(n)})\dd s(x).
\end{align*}
Hence, rearranging the terms in \eqref{eg8}, we obtain
\begin{align*}
	&\frac{\dd}{\dd t}\left(\frac12\int_\Omega|\nabla z^{(n)}|^2\dd x+\int_{\partial\Omega}\omega(x)\big(\hat\kappa_R(\theta^{(n)})-\theta^*K_R(\theta^{(n)})\big)\dd s(x)\right)\\
	&+\int_\Omega\mathcal{C}_V'(\theta^{(n)})\,\kappa(Q_R((\theta^{(n)})^+))\,|\theta^{(n)}_t|^2\dd x\\
	&= \int_\Omega\Bigg(\BB\nabla_s u^{(n)}_t:\nabla_s u^{(n)}_t+\|D_P[\nabla_s u^{(n)}]_t\|_*+\mu_R(p^{(n)})Q_R(|\nabla p^{(n)}|^2) \\
	& \qquad +(\chi^{(n)}+\rho^*(1-\chi^{(n)}))|D_0[p^{(n)}]_t|+\gamma_R(p^{(n)},\theta^{(n)},\dive u^{(n)})|\chi^{(n)}_t|^2\\
	& \qquad -\bigg(\frac{L}{\theta_c}\chi^{(n)}_t+\beta\dive u^{(n)}_t\bigg)Q_R((\theta^{(n)})^+)\Bigg)\kappa(Q_R((\theta^{(n)})^+))\theta^{(n)}_t\dd x-\int_{\partial\Omega}\omega(x)\theta^*_t K_R(\theta^{(n)})\dd s(x).
\end{align*}
We now integrate in time $\int_0^\tau \dd t$ for some $\tau \in [0,T_n)$. By Hypothesis \ref{h1}\,(viii) and (ix) it holds
\begin{align*}
	\int_\Omega|\nabla z^{(n)}|^2\dd x & \ge (\kappa^\flat)^2\int_\Omega|\nabla\theta^{(n)}|^2\dd x, \\[2 mm]
	\int_0^\tau\int_\Omega \mathcal{C}_V'(\theta^{(n)})\,\kappa(Q_R((\theta^{(n)})^+))\,|\theta^{(n)}_t|^2 \dd x\dd t &\ge c^\flat \kappa^\flat \int_0^\tau\io |\theta^{(n)}_t|^2 \dd x\dd t.
\end{align*}
Note also that $\hat\kappa_R(\theta) \ge \kappa^\flat\,\theta^2/2$ for all $\theta \in \real$. Hence, using Young's inequality as in \eqref{bou1} we obtain
\begin{align*}
	\int_{\partial\Omega}\omega(x)\big(\hat\kappa_R(\theta^{(n)})-\theta^*K_R(\theta^{(n)})\big)\dd s(x) &\ge \frac{\kappa^\flat}{2}\int_{\partial\Omega}\omega(x)|\theta^{(n)}|^2\dd s(x)-C_R\int_{\partial\Omega}\omega(x)|\theta^*\theta^{(n)}|^2\dd s(x)\\
	&\ge \frac{\kappa^\flat}{4}\int_{\partial\Omega}\omega(x)|\theta^{(n)}|^2\dd s(x)-C_R.
\end{align*}
Concerning the initial conditions, we employ Hypothesis \ref{h1}\,(iv) and (v). Thus, exploiting also \eqref{eneg2}, \eqref{enep2}, \eqref{gammaR_esti} and \eqref{esti0g}, we get
\begin{equation}\label{eg9}
	\begin{aligned}
	&\io |\nabla\theta^{(n)}|^2(x,\tau) \dd x + \int_{\partial\Omega}\omega(x)|\theta^{(n)}|^2(x,\tau)\dd s(x) + \int_0^\tau\io |\theta^{(n)}_t|^2(x,t)\dd t \\
	&\le C_R \, \bigg(1 + \int_0^\tau\io \left(|\nabla_s u^{(n)}_t|^2 + |p^{(n)}_t| + |p^{(n)}|^2 + |\dive u^{(n)}|^2\right)|\theta^{(n)}_t|\dd x\dd t \\
	& \quad + \int_0^\tau\int_{\partial\Omega} \omega(x)|\theta^*_t||\theta^{(n)}| \dd s(x)\dd t\bigg).	
	\end{aligned}
\end{equation}
We see that the approximate solution remains bounded in the maximal interval of existence $[0,T_n]$. Hence the solution exists globally, and for every $n \in \nat$ we have $T_n = T$.

We further need to estimate the terms $\nabla_s u^{(n)}_t$, $p^{(n)}$, $\dive u^{(n)}$ in the norm of $L^4(\Omega\times(0,T))$. Note that \eqref{esti3g}, \eqref{esti4g} entail $\nabla p^{(n)} \in L^\infty(0,T;L^2(\Omega;\real^3))$, $p^{(n)}_t \in L^2(\Omega\times(0,T))$ independently of $n$. Thus by the anisotropic embedding formulas (\cite[Theorem 10.2]{bin} on p.\,143 of the Russian version, see also \cite{adams,ee}) we deduce
\begin{equation}\label{eg10}
	\int_0^T\io |p^{(n)}|^4(x,t) \dd x\dd t \le C_R.
\end{equation}
Now, let us consider \eqref{e2g} rewritten in the form
\begin{equation}\label{rewr}
	\io P[\nabla_s u^{(n)}](x,t):\nabla_s\psi(x) \dd x + \io \BB\nabla_s u^{(n)}_t(x,t):\nabla_s\psi(x) \dd x = \io w^{(n)}(x,t)\,\dive\psi(x) \dd x,
\end{equation}
where
$$w^{(n)}(x,t) := p^{(n)}(\chi^{(n)}+\rho^*(1-\chi^{(n)}))+\beta(Q_R((\theta^{(n)})^+)-\theta_c)(x,t) + \widehat{G}(x,t)$$
according to Hypothesis \ref{h1}\,(ii). By the already mentioned $L^r$-regularity (with some $r \in [2,\infty)$) for elliptic systems in divergence form and by \eqref{ivP} we deduce, arguing as for \eqref{ex3},
\begin{equation}\label{nasut^p}
	\begin{aligned}
		\io |\nabla_s u^{(n)}_t|^r(x,t) \dd x &\le C \io |\nabla_s u^{(n)}|^r(x,0) \dd x + Ct^{r-1} \int_0^t\io |\nabla_s u_t^{(n)}|^r(x,\tau) \dd x \dd\tau \\
		& + C \io |w^{(n)}|^r(x,t) \dd x \quad \text{a.\,e.}
	\end{aligned}
\end{equation}
By \eqref{eg10} and Hypothesis \ref{h1}\,(ii) we see that
$$\int_0^\tau \io |w^{(n)}|^4(x,t) \dd x\dd t \le C_R \left(1 + \int_0^\tau \io |p^{(n)}|^4(x,t) \dd x\dd t\right) \le C_R.$$
Therefore, choosing $r=4$ in \eqref{nasut^p} and using also Hypothesis \ref{h1}\,(v), by Gr\"onwall's lemma we obtain
$$\int_0^\tau \io |\nabla_s u^{(n)}_t|^4(x,t) \dd x\dd t \le C_R,$$
from which we deduce a bound also for the term $\dive u^{(n)}$ since $\nabla_s u^{(n)}_t$ is dominant. Thus, coming back to \eqref{eg9} and using Hypothesis \ref{h1}\,(iii),\,(iv) and estimate \eqref{esti4g}, we finally obtain
\begin{align*}
	&\int_\Omega|\nabla\theta^{(n)}|^2(x,\tau)\dd x+\int_{\partial\Omega}\omega(x)|\theta^{(n)}|^2(x,\tau)\dd s(x)+\int_0^\tau\int_\Omega|\theta^{(n)}_t|^2(x,t)\dd x\dd t\cr
	&\le C_R\,\Bigg(1+\int_0^\tau\int_{\partial\Omega}\omega(x)|\theta^{(n)}|^2(x,t)\dd s(x)\dd t\Bigg).
\end{align*}
Applying Gr\"onwall's lemma and Poincar\'e's inequality \eqref{poin} we finally obtain the estimates
\begin{align}
	\supess_{\tau \in (0,T)} \left(\int_\Omega\left(|\theta^{(n)}|^2+|\nabla\theta^{(n)}|^2\right)(x,\tau)\dd x+\int_{\partial\Omega}\omega(x)|\theta^{(n)}|^2(x,\tau)\dd s(x)\right) &\le C_R, \label{esti6g} \\
	\int_0^T\int_\Omega |\theta^{(n)}_t|^2(x,t) \dd x\dd t &\le C_R. \label{esti7g}
\end{align}


\subsection{Limit as $n \to \infty$}\label{subsec_limn}

For the moment we keep the regularization parameters $\eta$ and $R$ fixed, and let $n \to \infty$ in \eqref{e1g}--\eqref{e4g}. From estimates \eqref{esti1g}, \eqref{esti2g}, \eqref{esti8g}, \eqref{esti3g}, \eqref{esti9g}, \eqref{esti4g}, \eqref{esti6g}, \eqref{esti7g}, we see that there exists a subsequence of $\{(p^{(n)},\theta^{(n)}) : n \in \nat\}$, which is again indexed by $n$, and functions $p,\theta$ such that
$$
\begin{array}{rcll}
	p^{(n)}_t \to p_t, \!\!\!\!\!&&\!\!\!\!\! \theta^{(n)}_t \to \theta_t &\quad\text{weakly in }L^2(\Omega \times (0,T)), \\
	\nabla\theta^{(n)} &\to& \nabla\theta &\quad\text{weakly-star in }L^\infty(0,T;L^2(\Omega;\real^3)), \\
	p^{(n)} &\to& p &\quad\text{strongly in }L^q(\Omega;C[0,T])\text{ for }q \in [1,6)\text{ and in }L^2(\partial\Omega \times (0,T)), \\
	\nabla p^{(n)} &\to& \nabla p &\quad\text{strongly in }L^2(\Omega\times(0,T);\real^3), \\
	\theta^{(n)} &\to& \theta &\quad\text{strongly in }L^2(\Omega \times (0,T))\text{ and in }L^2(\partial\Omega \times (0,T)),
\end{array}
$$
where the strong convergences are obtained by compact embedding, see \cite{bin}.
We also need strong convergence of the sequences $\{\nabla_s u^{(n)}\}$ and $\{\nabla_s u^{(n)}_t\}$ in order to pass to the limit in some nonlinear terms. Taking the difference of \eqref{rewr} for indices $n$ and $m$, and testing by $\psi = u^{(n)}_t-u^{(m)}_t$ we obtain, arguing as for Step 3 of the existence part,
\begin{equation}\label{nsutdiff}
	\begin{aligned}
		&\io |\nabla_s u^{(n)}_t-\nabla_s u^{(m)}_t|^2(x,\tau) \dd x \le C \io |w^{(n)}-w^{(m)}|^2(x,\tau) \dd x \\
		&\le C\left(\io \big(1+|p^{(n)}-p^{(m)}|^2+|p^{(n)}|^2|\chi^{(n)}-\chi^{(m)}|^2 +|\theta^{(n)}-\theta^{(m)}|^2\big)(x,\tau) \dd x\right)
	\end{aligned}
\end{equation}
a.\,e. in $(0,T)$. The $L^1$-Lipschitz continuity result for variational inequalities (see \cite[Theorem 1.12]{cmuc}) tells us that
\begin{equation}\label{lim1}
	|\chi^{(n)}-\chi^{(m)}|(x,\tau) \le \int_0^{\tau} |\chi^{(n)}_t-\chi^{(m)}_t|(x,t) \dd t \le 2 \int_0^{\tau} |F^{(n)}-F^{(m)}|(x,t) \dd t
\end{equation}
with the notation of \eqref{Fn}, where we have by virtue of \eqref{lip_G} for a.\,e. $x \in \Omega$ that
\begin{equation*}
	\begin{aligned}
		&\int_0^{\tau} |F^{(n)}-F^{(m)}|(x,t) \dd t \\
		&\le C_R \left(1 + \max_{s\in [0,\tau]}|p^{(n)}-p^{(m)}|(x,s) + \int_0^\tau \left(|\dive u^{(n)}-\dive u^{(m)}| + |\theta^{(n)}-\theta^{(m)}|\right)\!(x,t) \dd t\right).
	\end{aligned}
\end{equation*}
For $t \in [0,T]$ put
\begin{align*}
	U(t) &= \io |\nabla_s u^{(n)}-\nabla_s u^{(m)}|^2 (x,t) \dd x, \\
	W(t) &= \io |\nabla_s u^{(n)}_t-\nabla_s u^{(m)}_t|^2 (x,t) \dd x, \\
	\pi(t) &= 1 + \supess_{x \in \Omega} |p^{(n)}|^2(x,t), \\
	\Pi(t) &= \int_0^t \pi(s) \dd s, \\
	y(t) &= \io \left(1+\max_{s\in [0,t]}|p^{(n)}-p^{(m)}|^2(x,s)+|\theta^{(n)}-\theta^{(m)}|^2(x,t)\right) \dd x.
\end{align*}
Note that by \eqref{esti10g} $\Pi$ is bounded above independently of $n$. We further have
$$U(t) \le C\left(1 + \int_0^t W(s) \dd s\right),$$
and \eqref{nsutdiff} is of the form
$$W(\tau) \le C_R \left(y(\tau) + \pi(\tau) \int_0^\tau \left(y(t) + \int_0^t W(s) \dd s\right) \dd t\right).$$
Thus, from Fubini's theorem and Gr\"onwall's lemma we obtain
\begin{align*}
	\int_0^T W(\tau) \dd\tau &\le C_R \int_0^T \e^{C(\Pi(T)-\Pi(\tau))}\left(y(\tau) + \pi(\tau)\int_0^\tau y(t) \dd t\right) \dd \tau \\
	& \le C_R \left( \int_0^T \e^{C(\Pi(T)-\Pi(\tau))}y(\tau) \dd\tau + \int_0^T \left(\int_t^T \e^{C(\Pi(T)-\Pi(\tau))} \pi(\tau) \dd\tau\right) y(t) \dd t\right) \\
	&\le C_{R,\eta} \int_0^T y(\tau) \dd\tau,
\end{align*}
that is,
\begin{equation}\label{diff}
	\begin{aligned}
	&\io \supess_{\tau \in (0,T)}|\nabla_s u^{(n)}-\nabla_s u^{(m)}|^2(x,\tau) \dd x \le \int_0^T\io |\nabla_s u^{(n)}_t-\nabla_s u^{(m)}_t|^2(x,\tau) \dd x\dd\tau \\
	&\le C_{R,\eta} \left(1 + \io \left(\max_{s\in [0,T]}|p^{(n)}-p^{(m)}|^2(x,s)+\int_0^T|\theta^{(n)}-\theta^{(m)}|^2(x,\tau) \dd\tau\right) \dd x\right).
	\end{aligned}
\end{equation}
The sequences $\{p^{(n)}\}$ and $\{\theta^{(n)}\}$ are Cauchy in $L^2(\Omega;C[0,T])$ and $L^2(\Omega\times(0,T))$, respectively, hence $\{\nabla_s u^{(n)}\}$ and $\{\nabla_s u^{(n)}_t\}$ are also Cauchy sequences in $L^2(\Omega;C([0,T];\tens))$ and in $L^2(\Omega\times(0,T);\tens)$, respectively. Thus we conclude
$$
\begin{array}{rcll}
	\nabla_s u^{(n)} &\to& \nabla_s u &\quad \text{strongly in }L^2(\Omega;C([0,T];\tens)), \\
	\nabla_s u^{(n)}_t &\to& \nabla_s u_t &\quad \text{strongly in }L^2(\Omega\times(0,T);\tens).
\end{array}
$$
This is enough to pass to the limit in the nonlinearities by virtue of \cite[Theorem 12.10]{fuku}, and we obtain
$$
\begin{array}{rcll}
	f_R(p^{(n)}) \to f_R(p), \!\!\!\!&&\!\!\!\! \Phi_R(p^{(n)}) \to \Phi_R(p) & \ \text{strongly in }L^2(\Omega \times (0,T)), \\
	G_0[p^{(n)}] &\to& G_0[p] & \ \text{strongly in } L^2(\Omega;C[0,T]), \\
	G_0[p^{(n)}]_t &\to& G_0[p]_t & \ \text{weakly in }L^2(\Omega \times (0,T)), \\
	|D_0[p^{(n)}]_t| &\to& |D_0[p]_t| & \ \text{weakly in }L^2(\Omega \times (0,T)), \\
	U_0[p^{(n)}] &\to& U_0[p] & \ \text{strongly in } L^1(\Omega;C[0,T]), \\
	P[\nabla_s u^{(n)}] &\to& P[\nabla_s u] & \ \text{strongly in } L^2(\Omega;C([0,T];\tens)), \\
	\|D_P[\nabla_s u^{(n)}]_t\|_* &\to& \|D_P[\nabla_s u]_t\|_* & \ \text{weakly in }L^2(\Omega \times (0,T)), \\
	\BB\nabla_s u^{(n)}_t:\nabla_s u^{(n)}_t &\to& \BB\nabla_s u_t:\nabla_s u_t & \ \text{strongly in }L^1(\Omega \times (0,T)), \\
	\mathcal{C}_V(\theta^{(n)}) &\to& \mathcal{C}_V(\theta) & \ \text{strongly in }L^q(\Omega \times (0,T))\text{ for all }q \in \left[1,\frac{2}{1+\hat{b}}\right], \\
	&&\hspace{-4.9 cm}\left.
	\begin{array}{rcl}
		Q_R(|\nabla p^{(n)}|^2) &\to& Q_R(|\nabla p|^2) \\
		\mu_R(p^{(n)}) &\to& \mu_R(p) \\
		Q_R((\theta^{(n)})^+) &\to& Q_R(\theta^+) \\
		\displaystyle\frac{1}{\gamma_R(p^{(n)},\theta^{(n)},\dive u^{(n)})} &\to& \displaystyle\frac{1}{\gamma_R(p,\theta,\dive u)}
	\end{array}
	\right\rbrace & \ \text{strongly in }L^q(\Omega\times(0,T))\text{ for all }q\in[1,\infty),
\end{array}
$$
from which
$$
\begin{array}{rcll}
	f_R(p^{(n)})_t \to f_R(p)_t, \!\!\!\!&&\!\!\!\!
	\mathcal{C}_V(\theta^{(n)})_t \to \mathcal{C}_V(\theta)_t & \ \text{weakly in }L^2(\Omega \times (0,T)), \\	
	\gamma_R(p^{(n)},\theta^{(n)},\dive u^{(n)}) &\to& \gamma_R(p,\theta,\dive u) & \ \text{strongly in }L^q(\Omega\times(0,T))\text{ for all }q\in[1,\infty).
\end{array}
$$
The convergence of the Preisach hysteresis terms $G_0[p^{(n)}]$, $U_0[p^{(n)}]$ and $|D_0[p^{(n)}]_t|$ follow from \eqref{lip_G}, \eqref{cont_U0} and \eqref{en_id0}, respectively. The convergence of the plasticity terms $P[\nabla_s u^{(n)}]$ and $\|D_P[\nabla_s u^{(n)}]_t\|_*$ follows from \eqref{con1} and \eqref{enep1}, respectively. We now prove that the sequences $\{\chi^{(n)}\}$, $\{\chi^{(n)}_t\}$ converge strongly in appropriate function spaces. Inequality \eqref{lim1} yields
$$\io \sup_{\tau\in [0,T]}|\chi^{(n)}-\chi^{(m)}|(x,\tau) \dd x \le \int_0^{T}\io |\chi^{(n)}_t-\chi^{(m)}_t|(x,t) \dd x\dd t \le 2 \int_0^{T}\io |F^{(n)}-F^{(m)}|(x,t) \dd x\dd t$$
where, due to the above convergences, $F^{(n)} \to F$ strongly in $L^1(\Omega\times(0,T))$. Hence we conclude that $\{\chi^{(n)}(x,t)\}$ and $\{\chi^{(n)}_t(x,t)\}$ are Cauchy sequences in $L^1(\Omega;C[0,T])$ and in $L^1(\Omega\times(0,T))$, respectively. Moreover, since both $|\chi^{(n)}|$ and $|\chi^{(n)}_t|$ admit a uniform pointwise upper bound (see \eqref{esti0g}), we can use the Lebesgue dominated convergence theorem to conclude that
$$
\begin{array}{rcll}
	\chi^{(n)} &\to& \chi &\quad\text{strongly in }L^q(\Omega;C[0,T])\text{ for all }q\in[1,\infty), \\
	\chi^{(n)}_t &\to& \chi_t &\quad\text{strongly in }L^q(\Omega \times (0,T))\text{ for all }q\in[1,\infty).
\end{array}
$$
Passing to the limit as $n \to \infty$ in \eqref{e1g}--\eqref{e4g} we see that $(p,u,\theta,\chi)$ is a solution to \eqref{e1e}--\eqref{e4e}, \eqref{ini} with the regularity stated in Proposition~\ref{p1}.


\subsection{Limit as $\eta \to 0$}

Let us denote by $(v^{(\eta)},u^{(\eta)},z^{(\eta)},\chi^{(\eta)})$ the solution to \eqref{e1e}--\eqref{e4e}. Note that estimate \eqref{esti9g} is preserved in the limit as $n \to \infty$, hence the limit function $v^{(\eta)}$ satisfies the inequality
$$
\supess_{\tau \in (0,T)} \io |\Delta v^{(\eta)}|^2(x,\tau) \dd x \le \frac{C_R}{\eta}\,.
$$
We now choose $\phi \in W^{2,2}(\Omega)$ such that $\phi\big|_{\partial\Omega} = 0$, and integrate by parts in equation \eqref{e1e} to obtain 
\begin{equation}\label{e1e_parts}
	\begin{aligned}
		&\int_{\Omega}\big((\chi^{(\eta)}+\rho^*(1-\chi^{(\eta)}))(f_R(M_R^{-1}(v^{(\eta)}))+G_0[M_R^{-1}(v^{(\eta)})]+\dive u^{(\eta)})\big)_t\,\phi \dd x\\
		&+\int_{\Omega}\left(-\Delta v^{(\eta)}\,\phi + \eta\Delta v^{(\eta)}\Delta\phi\right) \dd x=0.
	\end{aligned}
\end{equation}
Introducing the new variable $\hat{v}^{(\eta)} = \Delta v^{(\eta)}$, we rewrite \eqref{e1e_parts} in the form
\begin{equation}\label{w_eq}
	\io \hat{v}^{(\eta)} \left(\phi-\eta\Delta\phi\right) \dd x = \io h^{(\eta)} \phi \dd x
\end{equation}
where
$$h^{(\eta)} = \big((\chi^{(\eta)}+\rho^*(1-\chi^{(\eta)}))(f_R(M_R^{-1}(v^{(\eta)}))+G_0[M_R^{-1}(v^{(\eta)})]+\dive u^{(\eta)})\big)_t\,.$$
Note that the term $G_0[p^{(\eta)}]_t$ is of order $p^{(\eta)}_t$ by \eqref{id_frak} and \eqref{G0}. Hence by estimates \eqref{esti0g}, \eqref{esti1g}, \eqref{esti2g}, \eqref{esti4g} we see that $h^{(\eta)} \in L^2(\Omega\times(0,T))$, and its $L^2$-norm is bounded independently of $\eta$.

Consider now the system $\{\hat e_k : k \in \nat\}$ of eigenfunctions of the negative Laplace operator with zero Dirichlet boundary conditions
$$-\Delta \hat e_k = \nu_k \hat e_k\,, \quad \hat e_k\big|_{\partial\Omega} = 0\,,\quad \io |\hat e_k(x)|^2\dd x = 1.$$
They form a complete orthonormal system in $L^2(\Omega)$ with $0 < \nu_1 \le \nu_2 \le \nu_3 \le \dots$. The functions $h^{(\eta)},\hat{v}^{(\eta)}$ admit the expansions
$$h^{(\eta)}(x,t) = \sum_{k=1}^\infty h^{(\eta)}_k(t)\,\hat e_k(x)\,, \quad \hat{v}^{(\eta)}(x,t) = \sum_{k=1}^\infty \hat{v}^{(\eta)}_k(t)\,\hat e_k(x)\,,$$
with coefficients $h^{(\eta)}_k : [0,T]\to\real$, $\hat{v}^{(\eta)}_k : [0,T] \to \real$. Choosing $\phi = \hat e_k$ in \eqref{w_eq} we obtain
$$\hat{v}^{(\eta)}_k(t) = \frac{h^{(\eta)}_k(t)}{1+\eta\nu_k}\,,$$
hence 
\begin{align}
	\int_0^T\io |\hat{v}^{(\eta)}(x,t)|^2\dd x\dd t &= \int_0^T\sum_{k=1}^\infty |\hat{v}^{(\eta)}_k(t)|^2\dd t \nonumber\\
	&\le \int_0^T\sum_{k=1}^\infty |h^{(\eta)}_k(t)|^2\dd t = \int_0^T\io |h^{(\eta)}(x,t)|^2\dd x\dd t \le C_R \label{limeta}
\end{align}
for some positive constant $C_R$ independent of $\eta$. The estimate \eqref{esti4g} is preserved in the limit $n \to \infty$. Thus we get for a subsequence $\eta \to 0$ that
$$
\begin{array}{rcll}
	\nabla v^{(\eta)} &\to& \nabla v &\quad\text{strongly in }L^2(\Omega\times(0,T);\real^3), \\
	\eta\,\Delta v^{(\eta)} &\to& 0 &\quad\text{strongly in } L^2(\Omega\times(0,T)),
\end{array}
$$
from which, by definition of $M_R$ in \eqref{MKR} and Hypothesis \ref{h1}\,(vii),
$$
\begin{array}{rcll}
	\nabla M_R^{-1}(v^{(\eta)}) &\to& \nabla M_R^{-1}(v) &\quad\text{strongly in }L^2(\Omega\times(0,T);\real^3), \\
	Q_R(|\nabla M_R^{-1}(v^{(\eta)})|^2) &\to& Q_R(|\nabla M_R^{-1}(v)|^2) &\quad\text{strongly in }L^q(\Omega\times(0,T))\text{ for all }q\in[1,\infty).
\end{array}
$$
Also the estimates \eqref{esti1g}, \eqref{esti2g}, \eqref{esti3g}, \eqref{esti6g}, \eqref{esti7g} are preserved when $n \to \infty$. Since they are independent of $\eta$, by letting $\eta \to 0$ we obtain for $v^{(\eta)}$ and $\theta^{(\eta)}$ the same convergences as before.

Note that as a side product, from \eqref{limeta} we get that the estimate
$$\int_0^T\io |\Delta M_R(p)|^2(x,t) \dd x\dd t \le C_R$$
holds also in the limit as $\eta \to 0$. Hence, arguing as for \eqref{esti8g} we get
$$\|M_R(p)\|^2_{L^2(0,T;W^{2,2}(\Omega))} \le C_R.$$
This, by Sobolev embedding, yields
\begin{equation}\label{l6}
	\int_0^T \left(\io |\nabla p|^6(x,t) \dd x\right)^{1/3} \dd t \le C_R,
\end{equation}
as well as
\begin{equation*}
	\int_0^T \supess_{x \in \Omega} |p^{(n)}|^2(x,t) \dd t \le C_R.
\end{equation*}
But then we can argue as in Subsection \ref{subsec_limn} and obtain an inequality similar to \eqref{diff}, but with a constant independent of $\eta$. This entails the strong convergence of the sequences $\nabla_s u^{(\eta)}$ and $\nabla_s u^{(\eta)}_t$. The rest of the convergence argument follows exactly as at the end of Subsection \ref{subsec_limn}, and this concludes the proof of Proposition \ref{p1}.

\section{Estimates independent of $R$}\label{estiR}

We now come back to our cut-off system \eqref{e1r}--\eqref{e4r}. We are going to derive a series of estimates independent of $R$. More precisely, after proving that the temperature stays away from zero, we will perform the energy estimate and the Dafermos estimate in order to gain some regularity for the temperature. Subsequently, a key-step will be the derivation of a bound for $p$ in an anisotropic Lebesgue space. Then an analogous estimate based on the particular structure of equation \eqref{e2r} is obtained for $\nabla_s u_t$. We finally show that this is sufficient for starting the Moser iteration and obtain an $L^\infty$ bound for $p$. After deriving some higher order estimates for the capillary pressure and for the temperature, we will be ready to let $R$ tend to $\infty$ in \eqref{e1r}--\eqref{e4r}.


\subsection{Positivity of the temperature}

For every nonnegative test function $\zeta \in X$ we have, by virtue of \eqref{e3r},
\begin{align*}
	&\int_\Omega\left(\mathcal{C}_V(\theta)_t\zeta+\kappa(Q_R(\theta^+))\nabla\theta\cdot\nabla\zeta\right)\dd x+\int_{\partial\Omega}\omega(x)(\theta-\theta^*)\zeta\dd s(x) \cr
	&= \int_\Omega\bigg(\BB\nabla_s u_t:\nabla_s u_t+\|D_P[\nabla_s u]_t\|_*+\mu_R(p)Q_R(|\nabla p|^2)+(\chi+\rho^*(1-\chi))|D_0[p]_t| \cr
	& \quad +\gamma_R(p,\theta,\dive u)\chi_t^2-\left(\frac{L}{\theta_c}\chi_t+\beta\dive u_t\right)Q_R(\theta^+)\bigg)\zeta\dd x \cr
	&\ge \int_\Omega\left(\frac{B^\flat}{3}|\dive u_t|^2+\gamma^\flat\chi_t^2-\bigg(\frac{L}{\theta_c}\chi_t+\beta\dive u_t\bigg)Q_R(\theta^+)\right)\zeta\dd x,
\end{align*}
where in the last line we used Hypothesis \ref{h1}\,(i) together with inequality \eqref{div}, and also estimate \eqref{gammaR_esti}. Then, by Young's inequality,
$$\int_\Omega\left(\mathcal{C}_V(\theta)_t\zeta+\kappa(Q_R(\theta^+))\nabla\theta\cdot\nabla\zeta\right)\dd x+\int_{\partial\Omega}\omega(x)(\theta-\theta^*)\zeta\dd s(x) \ge -C\int_\Omega(Q_R(\theta^+))^2\zeta\dd x$$
with a constant $C$ depending on $L,\,\theta_c,\,\beta,\,\BB,\,\gamma^\flat$. Let now $\varphi(t)$ be the solution of the ODE
$$\frac{\dd}{\dd t}\mathcal{C}_V(\varphi(t)) + C\varphi^2(t) = 0, \quad \varphi(0)=\bar\theta$$
with $\bar\theta$ from Hypothesis \ref{h1}. Then $\varphi$ is nondecreasing and positive. Taking into account the fact that $\mathcal{C}_V(\varphi)_t = -C\varphi^2$ and $\nabla\varphi=0$, for every nonnegative test function $\zeta\in X$ we have in particular
\begin{align*}
	&\int_\Omega\left(\big(\mathcal{C}_V(\varphi)-\mathcal{C}_V(\theta)\big)_t\,\zeta+\kappa(Q_R(\theta^+))\nabla(\varphi-\theta)\cdot\nabla\zeta\right)\dd x+\int_{\partial\Omega}\omega(x)(\theta-\theta^*)\zeta\dd s(x)\\
	&\le C\int_\Omega\left((Q_R(\theta^+))^2-\varphi^2\right)\zeta\dd x.
\end{align*}
Consider now the following regularization of the Heaviside function
$$
H_{\ve}(z) =
\left\lbrace
\begin{aligned}
	0 &\quad\text{for }z \le 0,\\
	\frac{z}{\ve} &\quad\text{for }0<z\le{\ve},\\
	1 &\quad\text{for }z>{\ve},
\end{aligned}
\right.
$$
for $\ve > 0$, and set $\zeta(x,t)=H_{\ve}(\varphi(t)-\theta(x,t))$ which is an admissible test function. This yields
$$\io \big(\mathcal{C}_V(\varphi)-\mathcal{C}_V(\theta)\big)_t\,H_{\ve}(\varphi-\theta) \dd x \le 0.$$
By the Lebesgue Dominated Convergence Theorem we can pass to the limit in the above inequality for $\ve \to 0$, getting
$$\io \big(\mathcal{C}_V(\varphi)-\mathcal{C}_V(\theta)\big)_t\,H(\varphi-\theta) \dd x \le 0,$$
that is, by the monotonicity of $\mathcal{C}_V$,
$$\frac{\dd}{\dd t} \int_\Omega \left(\mathcal{C}_V(\varphi)-\mathcal{C}_V(\theta)\right)^+\dd x\le 0, \qquad (\mathcal{C}_V(\varphi)-\mathcal{C}_V(\theta))^+(x,0)=0$$
which implies $(\mathcal{C}_V(\varphi)-\mathcal{C}_V(\theta))^+ \equiv 0$. Owing again to the monotonicity of $\mathcal{C}_V$ and $\varphi$, we conclude that, independently of $R$,
\begin{equation}\label{pos}
	\theta(x,t) \ge \varphi(t) \ge \varphi(T) =: \theta_T>0 \quad \text{for all }x\text{ and }t.
\end{equation}
We now pass to a series of estimates independent of $R$.


\subsection{Energy estimate}

Since we proved that the temperature stays positive, from now on we will write $Q_R(\theta^+)=Q_R(\theta)$.

We test \eqref{e1r} by $\phi=p$, \eqref{e2r} by $\psi=u_t$ and \eqref{e3r} by $\zeta=1$. Summing up the three resulting equations we obtain
\begin{align*}
	&\int_\Omega \big((\chi+\rho^*(1-\chi))(f_R(p)+G_0[p]+\dive u)\big)_t\,p \dd x + \int_\Omega \mu_R(p)|\nabla p|^2 \dd x \\
	& + \int_\Omega (P[\nabla_s u]+\BB\nabla_s u_t):\nabla_s u_t \dd x - \int_\Omega \Big(p\,(\chi+\rho^*(1-\chi))+\beta(Q_R(\theta)-\theta_c)\Big)\,\dive u_t \dd x \\
	& + \int_\Omega \bigg(\mathcal{C}_V(\theta)_t-\BB\nabla_s u_t:\nabla_s u_t-\|D_P[\nabla_s u]_t\|_*-\mu_R(p)Q_R(|\nabla p|^2) \\
	&-(\chi+\rho^*(1-\chi))|D_0[p]_t|-\gamma_R(p,\theta,\dive u)\chi_t^2\bigg) \dd x + \int_\Omega \bigg(\frac{L}{\theta_c}\chi_t+\beta\dive u_t\bigg)Q_R(\theta) \dd x \\
	&= \int_{\partial\Omega} \alpha(x)(p^*-p)p \dd s(x) + \int_\Omega g\cdot u_t \dd x + \int_{\partial\Omega} \omega(x)(\theta^*-\theta) \dd s(x).
\end{align*}
Note that some of the terms cancel out. Moreover, recalling the notation introduced in \eqref{R} and the energy balance \eqref{en_id0}, the identities
\begin{align}
	&
	\begin{aligned}
		&\int_\Omega\big((\chi+\rho^*(1-\chi))f_R(p)\big)_t\,p\dd x \\
		&= \frac{\dd}{\dd t}\int_\Omega(\chi+\rho^*(1-\chi))V_R(p)\dd x+\int_\Omega(1-\rho^*)\chi_t\Phi_R(p)\dd x,
	\end{aligned}
	\label{il}
	\\[2 mm]
	&
	\begin{aligned}
		&\int_\Omega\big((\chi+\rho^*(1-\chi))G_0[p]\big)_t\,p\dd x - \io (\chi+\rho^*(1-\chi))|D_0[p]_t| \dd x \\
		& = \frac{\dd}{\dd t}\int_\Omega(\chi+\rho^*(1-\chi))U_0(p)\dd x + \io (1-\rho^*)\chi_t\left(p\,G_0[p]-U_0[p]\right) \dd x
	\end{aligned}
	\label{lo}
\end{align}
hold true. Hence we obtain, using also \eqref{enep1} and \eqref{e4r},
\begin{equation}\label{ee1}
	\begin{aligned}
		&\frac{\dd}{\dd t}\int_\Omega \Big(\mathcal{C}_V(\theta)+L\chi+\beta\theta_c\,\dive u+(\chi+\rho^*(1-\chi))(V_R(p)+U_0[p])+U_P[\nabla_s u]\Big) \dd x \\
		& + \int_\Omega \mu_R(p)\Big(|\nabla p|^2-Q_R(|\nabla p|^2)\Big) \dd x + \int_{\partial\Omega} \Big(\alpha(x)(p-p^*)p+\omega(x)(\theta-\theta^*)\Big) \dd s(x) \\
		& = \frac{\dd}{\dd t} \int_\Omega g\cdot u \dd x - \int_\Omega g_t\cdot u \dd x.
	\end{aligned}
\end{equation}
We now integrate in time $\int_0^\tau \dd t$. On the left-hand side Young's inequality, \eqref{enep2} and \eqref{div} entail
\begin{equation}\label{d4}
	U_P[\nabla_s u] + \beta\theta_c\,\dive u \ge \frac{A^\flat}{8}|\nabla_s u|^2 - C.
\end{equation}
By the definition of $Q_R$ in \eqref{QR}, it holds $|\nabla p|^2 \ge Q_R(|\nabla p|^2)$. The boundary term is such that
\begin{equation}\label{ee2}
\int_0^\tau\int_{\partial\Omega} \Big(\alpha(x)(p-p^*)p+\omega(x)(\theta-\theta^*)\Big) \dd s(x) \dd t \ge \int_0^\tau\int_{\partial\Omega} \left(\alpha(x)\frac{p^2}{2}+\omega(x)\theta\right)\dd s(x)\dd t - C
\end{equation}
thanks to Young's inequality and Hypothesis \ref{h1}\,(iii) and (iv). Concerning the right-hand side of \eqref{ee1}, the time integration gives
\begin{align*}
	&\int_\Omega (g \cdot u)(x,\tau) \dd x - \int_\Omega g(x,0) \cdot u^0(x) \dd x - \int_0^\tau\io (g_t \cdot u)(x,t) \dd x\dd t,
\end{align*}
where the term containing the initial conditions is controlled by using H\"older's inequality and observing that
$$\io |g|^2(x,0) \dd x = \io |g|^2(x,\tau) \dd x - 2\int_0^\tau\io (g \cdot g_t)(x,t) \dd x\dd t.$$
Hence by Young's inequality and Hypothesis \ref{h1}\,(ii), (v) we deduce
\begin{align*}
	&\int_\Omega (g \cdot u)(x,\tau) \dd x - \int_\Omega g(x,0) \cdot u^0(x) \dd x - \int_0^\tau\io (g_t \cdot u)(x,t) \dd x\dd t \\
	&\le \frac{A^\flat}{16}\io |\nabla_s u|^2(x,\tau) \dd x + C \left(1 + \int_0^\tau\io |\nabla_s u|^2(x,t) \dd x\dd t\right),
\end{align*}
where we used also Korn's inequality \eqref{korn}. The first term in the last line is absorbed by \eqref{d4}. Finally, the initial conditions are kept under control thanks to \eqref{eneg2}, \eqref{enep}, \eqref{fR_esti} and Hypothesis \ref{h1}\,(v). Hence what we eventually get is
\begin{align*}
	&\int_\Omega \left(\mathcal{C}_V(\theta)+V_R(p)+|\nabla_s u|^2\right)(x,\tau) \dd x + \int_0^\tau\int_{\partial\Omega} \Big(\alpha(x)p^2+\omega(x)\theta\Big)(x,t) \dd s(x)\dd t \\
	&\le C \left( 1 + \int_0^\tau\int_\Omega |\nabla_s u|^2(x,t) \dd x\dd t \right),
\end{align*}
and applying Gr\"onwall's lemma we finally obtain the estimates
\begin{align}
	\supess_{\tau \in (0,T)} \io \left(\mathcal{C}_V(\theta)+V_R(p)+|\nabla_s u|^2\right)(x,\tau) \dd x &\le C, \label{esti_en1} \\
	\int_0^T\int_{\partial\Omega} \Big(\alpha(x)p^2+\omega(x)\,\theta\Big)(x,t) \dd s(x)\dd t &\le C. \label{esti_en2}
\end{align}
Estimate \eqref{esti_en1} also gives
\begin{equation}\label{esti_en3}
	\supess_{\tau \in (0,T)} \io |\theta|^{1+b}(x,\tau) \dd x \le C,
\end{equation}
where $b$ is from Hypothesis \ref{h1}\,(viii).


\subsection{Dafermos estimate}

We set $\hat\theta := Q_R(\theta)$ and test \eqref{e3r} by $\zeta = -\hat\theta^{-a}$, with $a$ from Hypothesis \ref{h1}. This yields the identity
\begin{equation}\label{de1}
	\begin{aligned}
		&\int_{\Omega}\bigg(\mathcal{C}_V(\theta)_t-\BB\nabla_s u_t:\nabla_s u_t-\|D_P[\nabla_s u]_t\|_*-\mu_R(p)Q_R(|\nabla p|^2)-\dch|D_0[p]_t| \\
		& \ -\gamma_R(p,\theta,\dive u)\chi_t^2+\left(\frac{L}{\theta_c}\chi_t+\beta\dive u_t\right)\hat\theta\bigg)(-\hat\theta^{-a}) \dd x + \int_{\Omega}\kappa(\hat\theta)\nabla\theta\cdot\nabla(-\hat\theta^{-a}) \dd x \\
		& = \int_{\partial\Omega}\omega(x)(\theta^*-\theta)\,(-\hat\theta^{-a}) \dd s(x).
	\end{aligned}
\end{equation}
It holds
$$\io \mathcal{C}_V(\theta)_t\,(-\hat\theta^{-a}) \dd x = -\frac{\dd}{\dd t}\io F_a(\theta) \dd x$$
where
$$F_a(\theta) := \int_0^\theta \frac{c_V(s)}{(Q_R(s))^a} \dd s,$$
and by Hypothesis \ref{h1}\,(ix) also
$$\io \kappa(\hat{\theta})\nabla\theta\cdot\nabla(-\hat{\theta}^{-a}) \dd x = \io \kappa(\hat{\theta})\,a\,\hat{\theta}^{-a-1}\nabla\theta\cdot\nabla\hat{\theta} \dd x \ge a\kappa^\flat \io |\nabla\hat{\theta}|^2 \dd x.$$
Hence from \eqref{de1} we get, using also Hypothesis \ref{h1}\,(i) and inequalities \eqref{div}, \eqref{gammaR_esti},
\begin{equation}\label{de2}
	\begin{aligned}
		&\io \bigg(\frac{B^\flat}{3}|\dive u_t|^2+\gamma^\flat\chi_t^2\bigg)\hat\theta^{-a} \dd x + a\kappa^\flat\int_{\Omega}|\nabla\hat\theta|^2 \dd x \\
		&\le \io \left(\frac{L}{\theta_c}\chi_t+\beta\dive u_t\right)\hat\theta^{1-a} \dd x + \int_{\partial\Omega}\omega(x)(\theta-\theta^*)\,\hat\theta^{-a} \dd s(x) + \frac{\dd}{\dd t} \io F_a(\theta) \dd x,
	\end{aligned}
\end{equation}
where we neglected some positive terms on the left-hand side. Young's inequality yields
$$\left(\frac{L}{\theta_c}\chi_t+\beta\dive u_t\right)\hat\theta^{1-a} \le \left(\frac{\gamma^\flat}{2}\chi_t^2 + \frac{B^\flat}{4}|\dive u_t|^2\right)\hat\theta^{-a} + C\hat\theta^{2-a},$$
with a constant $C$ depending only on $L,\,\theta_c,\,\beta,\,\BB,\,\gamma^\flat$, whereas the boundary term is such that
$$\int_{\partial\Omega}\omega(x)(\theta-\theta^*)\,\hat\theta^{-a}\dd s(x) \le C\left(1 + \int_{\partial\Omega}\omega(x)\theta \dd s(x)\right)$$
by \eqref{pos} and Hypothesis \ref{h1}\,(iii),\,(iv). Note also that $F_a(\theta) \le \mathcal{C}_V(\theta)$ for all $\theta \ge 0$. Thus, integrating \eqref{de2} in time $\int_0^\tau \dd t$ for some $\tau \in [0,T]$ and neglecting some other positive terms on the left-hand side we obtain
\begin{equation}\label{daf1}
	\int_0^\tau\io |\nabla\hat\theta|^2(x,t) \dd x \dd t \le C \left(1 + \int_0^\tau\io \hat\theta^{2-a}(x,t) \dd x \dd t\right)
\end{equation}
thanks to estimates \eqref{esti_en1}, \eqref{esti_en2}. Now, owing to estimate \eqref{esti_en3}, we can apply the Gagliardo-Nirenberg inequality \eqref{gn} with the choices $s=1+b$, $r=2$ and $N=3$ obtaining, for $t\in(0,T)$,
$$|\hat{\theta}(t)|_q \le C \left( 1 + |\nabla\hat{\theta}(t)|_2^\delta \right)$$
with $\delta=\frac{6(q-1-b)}{(5-b)q}$ and for every $1+b<q<6$. In particular, since $\delta\cdot\frac{(5-b)q}{3(q-1-b)}=2$, this and \eqref{daf1} yield
\begin{equation}\label{daf2}
	\int_0^T\left(|\hat{\theta}(t)|_q^q\right)^{(5-b)/3(q-1-b)}\dd t \le C \left(1 + \int_0^T|\nabla\hat{\theta}(t)|_2^2\dd t\right) \le C \left(1 + \int_0^T|\hat{\theta}(t)|_{2-a}^{2-a}\dd t\right).
\end{equation}
Let us now choose $q=2-a$, which is admissible in the sense that $1+b<2-a$ thanks to Hypothesis \ref{h1}. Since $\frac{5-b}{3(1-a-b)}>1$, we can apply Young's inequality on the right-hand side getting
$$\int_0^T |\hat{\theta}(t)|_{2-a}^{2-a} \dd t \le C.$$
Substituting in \eqref{daf2} entails
\begin{equation}\label{esti_daf1}
	\int_0^T\io |\nabla\hat\theta|^2(x,t) \dd x \dd t \le C.
\end{equation}
Coming back to \eqref{daf2} again and choosing $q=8/3+2b/3$, we also get
\begin{equation}\label{esti_daf2}
	\int_0^T\io \hat\theta^{8/3+2b/3}(x,t) \dd x \dd t \le C.
\end{equation}


\subsection{Mechanical energy estimate}\label{subsec_mechen}

In order to estimate the capillary pressure in a suitable anisotropic Lebesgue space, we first need to find a bound for $\dive u_t$ in $L^2(\Omega\times(0,T))$, independently of $R$. To this purpose, we test \eqref{e1r} by $\phi = p$, \eqref{e2r} by $\psi = u_t$ and sum up to obtain, with the notation of the previous subsection,
\begin{align*}
	&\io \big((\chi+\rho^*(1-\chi))(f_R(p)+G_0[p]+\dive u)\big)_t\,p \dd x + \io \mu_R(p)|\nabla p|^2 \dd x \\
	& + \io (P[\nabla_s u]+\BB\nabla_s u_t):\nabla_s u_t \dd x - \io \big(p(\chi+\rho^*(1-\chi))+\beta(\hat\theta-\theta_c)\big)\,\dive u_t \dd x \\
	&= \int_{\partial\Omega} \alpha(x)(p^*-p)\,p \dd s(x) + \io g\cdot u_t \dd x.
\end{align*}
Note that some terms cancel out. Owing to \eqref{il} and exploiting also the energy identity \eqref{enep1}, what we eventually get is
\begin{align*}
	&\frac{\dd}{\dd t} \io \Big((\chi+\rho^*(1-\chi))(V_R(p)+U_0[p])+U_P[\nabla_s u]\Big) \dd x + \io \mu_R(p)|\nabla p|^2 \dd x \\
	& + \io \BB\nabla_s u_t:\nabla_s u_t \dd x + \io (1-\rho^*)\chi_t\,\big(\Phi_R(p)+p\,G_0[p]+p\,\dive u-U_0[p]\big) \dd x \\
	& + \int_{\partial\Omega} \alpha(x)(p-p^*)\,p \dd s(x) \le \io \beta(\hat\theta-\theta_c)\,\dive u_t \dd x + \io g\cdot u_t \dd x.
\end{align*}
Now, \eqref{e4r} yields
\begin{align}
	&(1-\rho^*)\chi_t\,\big(\Phi_R(p)+p\,G_0[p]+p\,\dive u-U_0[p]\big) = \gamma_R(p,\theta,\dive u)\chi_t^2 - L\chi_t\left(\frac{\hat\theta}{\theta_c}-1\right) \nonumber\\
	&= \gamma_R(p,\theta,\dive u)\chi_t^2 - \sqrt{\gamma_R(p,\theta,\dive u)}\chi_t\,\frac{L}{\sqrt{\gamma_R(p,\theta,\dive u)}}\left(\frac{\hat\theta}{\theta_c}-1\right) \nonumber\\
	& \ge \frac12 \gamma_R(p,\theta,\dive u)\chi_t^2 - C(1+\hat\theta) \label{me1}
\end{align}
where in the last line we used Young's inequality and \eqref{gammaR_esti}, and where the constant $C$ is independent of $R$. Moreover, from the pointwise inequality \eqref{div} and arguing as for \eqref{eg3} we get
\begin{align*}
	\io \beta(\hat\theta-\theta_c)\,\dive u_t \dd x & \le C \left(1 + \io \hat\theta^2 \dd x\right) + \frac{B^\flat}{4}\io |\nabla_s u_t|^2 \dd x, \\
	\io g\cdot u_t \dd x &\le C + \frac{B^\flat}{4} \io |\nabla_s u_t|^2 \dd x.
\end{align*}
Hence we obtain, exploiting also Hypothesis \ref{h1}\,(i) to absorb the terms coming from the two estimates above,
\begin{align*}
	&\frac{\dd}{\dd t} \io \Big((\chi+\rho^*(1-\chi))(V_R(p)+U_0[p])+U_P[\nabla_s u]\Big) \dd x + \mu^\flat\io |\nabla p|^2 \dd x + \frac{B^\flat}{2}\io |\nabla_s u_t|^2 \dd x \\
	& + \frac12 \io \gamma_R(p,\theta,\dive u)\chi_t^2 \dd x + \frac12 \int_{\partial\Omega} \alpha(x)p^2 \dd s(x) \le C \left(1 + \io \hat\theta^2 \dd x\right)
\end{align*}
where the boundary term was handled as in \eqref{ee2}. We now integrate in time $\int_0^\tau \dd t$ for some $\tau \in [0,T]$. The right-hand side is bounded thanks to estimate \eqref{esti_daf2}, whereas the initial conditions are kept under control thanks to \eqref{eneg2}, \eqref{vari}, \eqref{enep}, \eqref{fR_esti} and Hypothesis \ref{h1}\,(v). Hence, neglecting some already estimated positive terms, we finally obtain
\begin{equation}\label{esti5}
	\int_0^T\io \left(|\nabla p|^2+|\nabla_s u_t|^2\right)(x,t) \dd x\dd t \le C
\end{equation}
independently of $R$. This, together with \eqref{esti_en2} and Poincar\'e's inequality \eqref{poin}, yields
\begin{equation}\label{esti6}
	\|p\|_{L^2(0,T;W^{1,2}(\Omega))}^2 \le C.
\end{equation}


\subsection{Estimate for the capillary pressure}\label{subsec_estp}

We choose an even function $\lambda : \real \to (0,\infty)$ such that $\lambda'(p) \ge 0$ for $p>0$ and $p\lambda(p) \in X$. Then we test \eqref{e1r} by $\phi = p\lambda(p)$. We obtain
\begin{equation}\label{cp0}
	\begin{aligned}
		&\int_{\Omega}\big((\chi+\rho^*(1-\chi))(f_R(p)+G_0[p]+\dive u)\big)_t\,p\lambda(p) \dd x + \int_{\Omega}\mu_R(p)(\lambda(p)+p\lambda'(p))|\nabla p|^2 \dd x \\
		&= \int_{\partial\Omega}\alpha(x)(p^*-p)\,p\lambda(p) \dd s(x).
	\end{aligned}
\end{equation}
The term under the time derivative has the form
\begin{equation}\label{cp1}
	\begin{aligned}
		&\io \big((\chi+\rho^*(1-\chi))(f_R(p)+G_0[p]+\dive u)\big)_t\,p\lambda(p) \dd x \\
		&= \io (1-\rho^*)\chi_t\left(f_R(p)+G_0[p]+\dive u\right)p\lambda(p) \dd x + \io (\chi+\rho^*(1-\chi))f_R'(p)p_t\,p\lambda(p) \dd x \\
		& + \io (\chi+\rho^*(1-\chi))\,G_0[p]_t\,p\lambda(p) \dd x + \io (\chi+\rho^*(1-\chi))\,\dive u_t\,p\lambda(p) \dd x.
	\end{aligned}
\end{equation}
We now define
$$V_{\lambda,R}(p) := \int_0^p f_R'(z)z\lambda(z) \dd z$$
so that
$$\io (\chi+\rho^*(1-\chi))f_R'(p)p_t\,p\lambda(p) \dd x = \frac{\dd}{\dd t} \io (\chi+\rho^*(1-\chi))V_{\lambda,R}(p) \dd x - \io (1-\rho^*)\chi_t V_{\lambda,R}(p) \dd x,$$
and introduce the modified Preisach potential as a counterpart to \eqref{Uh}
$$U_\lambda[p] := \int_0^\infty\int_0^{\play_r[p]} v\lambda(v)\,\psi(r,v)\dd v\dd r$$
which satisfies
$$G_0[p]_t\,p\lambda(p)-U_\lambda[p]_t \ge 0 \ \text{ a.\,e.}$$
according to \eqref{Uh_ineq}. Note that $V_{\lambda,R}(p) > 0$ and $U_\lambda[p] \ge 0$ for all $p \ne 0$. Then \eqref{cp1} can be rewritten as
\begin{equation}\label{cp2}
	\begin{aligned}
		&\io \big((\chi+\rho^*(1-\chi))(f_R(p)+G_0[p]+\dive u)\big)_t\,p\lambda(p) \dd x \\
		& \ge \frac{\dd}{\dd t} \io (\chi+\rho^*(1-\chi))\left(V_{\lambda,R}(p)+U_\lambda[p]\right) \dd x + \io (\chi+\rho^*(1-\chi))\,\dive u_t\,p\lambda(p) \dd x \\
		& + \io (1-\rho^*)\chi_t\,\Big(\left(pf_R(p)+p\,G_0[p]+p\,\dive u\right)\lambda(p)-V_{\lambda,R}(p)-U_\lambda[p]\Big) \dd x.
	\end{aligned}
\end{equation}
Defining
$$\Psi_{\lambda,R}(p) := V_R(p)\lambda(p) - V_{\lambda,R}(p),$$
we see that from \eqref{R} and \eqref{e4r} it holds
\begin{align*}
	&(1-\rho^*)\chi_t\,\Big(\left(pf_R(p)+p\,G_0[p]+p\,\dive u\right)\lambda(p)-V_{\lambda,R}(p)-U_\lambda[p]\Big) \\
	&= (1-\rho^*)\chi_t\,\Big(\left(\Phi_R(p)+p\,G_0[p]-U_0[p]+p\,\dive u\right)\lambda(p)+\Psi_{\lambda,R}(p)+U_0[p]\lambda(p)-U_\lambda[p]\Big) \\
	&= \left(\gamma_R(p,\theta,\dive u)\chi_t^2 - L\chi_t\left(\frac{\hat\theta}{\theta_c}-1\right)\right)\lambda(p) + (1-\rho^*)\chi_t\left(\Psi_{\lambda,R}(p)+U_0[p]\lambda(p)-U_\lambda[p]\right).
\end{align*}
Now, using Young's inequality as in \eqref{me1} we obtain
$$\left(\gamma_R(p,\theta,\dive u)\chi_t^2 - L\chi_t\left(\frac{\hat\theta}{\theta_c}-1\right)\right)\lambda(p) \ge \frac12 \gamma_R(p,\theta,\dive u)\chi_t^2\,\lambda(p) - C(1+\hat\theta)\,\lambda(p),$$
and similarly
\begin{align*}
	&\big|(1-\rho^*)\chi_t\left(\Psi_{\lambda,R}(p)+U_0[p]\lambda(p)-U_\lambda[p]\right)\!\big| \\
	&\le \frac14 \gamma_R(p,\theta,\dive u)\chi_t^2\,\lambda(p) + C\,\frac{\left(\Psi_{\lambda,R}(p)+U_0[p]\lambda(p)-U_\lambda[p]\right)^2}{\gamma_R(p,\theta,\dive u)\,\lambda(p)}\,,
\end{align*}
so that \eqref{cp2} entails
\begin{equation}\label{cp4}
	\begin{aligned}
		&\io \big((\chi+\rho^*(1-\chi))(f_R(p)+G_0[p]+\dive u)\big)_t\,p\lambda(p) \dd x \\
		& \ge \frac{\dd}{\dd t} \io (\chi+\rho^*(1-\chi))\left(V_{\lambda,R}(p)+U_\lambda[p]\right) \dd x + \io (\chi+\rho^*(1-\chi))\,\dive u_t\,p\lambda(p) \dd x \\
		& + \frac14 \io \gamma_R(p,\theta,\dive u)\chi_t^2\,\lambda(p) \dd x - C \io \left((1+\hat\theta)\,\lambda(p) + \frac{\left(\Psi_{\lambda,R}(p)+U_0[p]\lambda(p)-U_\lambda[p]\right)^2}{\gamma_R(p,\theta,\dive u)\,\lambda(p)}\right) \dd x.
	\end{aligned}
\end{equation}
Note that
$$V_R(p) \le V(p)+\frac{f^\sharp}{2}(p^2-R^2)^+$$
for all $p \in \real$, hence
\begin{align*}
	\frac{\left(\Psi_{\lambda,R}(p)+U_0[p]\lambda(p)-U_\lambda[p]\right)^2}{\left(1+p^2\right)\gamma_R(p,\theta,\dive u)\,\lambda^2(p)} &\le \frac{\left(V_R(p)+U_0[p]\right)^2}{(1+p^2)\,\gamma_R(p,\theta,\dive u)} \\
	&\le \frac{\left(V(p)+\frac{f^\sharp}{2}(p^2-R^2)^++U_0[p]\right)^2}{\gamma^\flat(1+p^2)\left(1+(p^2-R^2)^+\right)} \le C
\end{align*}
independently of $R$. From \eqref{cp4} we conclude
\begin{align*}
	&\io \big((\chi+\rho^*(1-\chi))(f_R(p)+G_0[p]+\dive u)\big)_t\,p\lambda(p) \dd x \\
	& \ge \frac{\dd}{\dd t} \io (\chi+\rho^*(1-\chi))\left(V_{\lambda,R}(p)+U_\lambda[p]\right) \dd x + \frac14 \io \gamma_R(p,\theta,\dive u)\chi_t^2\,\lambda(p) \dd x \\
	& - C \io \left(1+|p||\dive u_t|+\hat\theta+p^2\right)\lambda(p) \dd x,
\end{align*}
so that \eqref{cp0} and a time integration $\int_0^\tau \dd t$ for some $\tau \in [0,T]$ yields
\begin{equation}\label{cp3}
	\begin{aligned}
		&\io (\chi+\rho^*(1-\chi))\left(V_{\lambda,R}(p)+U_\lambda[p]\right)\!(x,\tau) \dd x \\
		& + \int_0^\tau\io \mu_R(p)(\lambda(p)+p\lambda'(p))|\nabla p|^2(x,t) \dd x\dd t + \int_0^\tau\int_{\partial\Omega} \alpha(x)(p-p^*)\,p\lambda(p)(x,t) \dd s(x)\dd t \\
		&\le C \int_0^\tau\io \left(1+|p||\dive u_t|+\hat\theta+p^2\right)\lambda(p)(x,t) \dd x\dd t + \io \left(V_{\lambda,R}(p)+U_\lambda[p] \right)\!(x,0) \dd x.
	\end{aligned}
\end{equation}
Now that we got rid of $\chi_t$ and derived a manageable estimate, we choose $\lambda(p) = |p|^{2k}$ with $k \ge \nu/2$ which will be specified later. Here $\nu$ is as in Hypothesis \ref{h1}. Note that this is an admissible choice, that is, $p\lambda(p) \in X$. Indeed, by estimates \eqref{esti4g}, \eqref{l6} and by the anisotropic
embedding formulas, see \cite{bin}, we have
\begin{equation}\label{ani1}
	p \in L^q(0,T;C(\bar\Omega))
\end{equation}
for any $q \in [1,4)$. The bound also depends on $R$, but for a fixed $R$ and each $k>0$ the function $p|p|^{2k}(\cdot, t)$ belongs to $X$ for a.\,e. $t \in (0,T)$.

With this choice \eqref{cp3} takes the form
\begin{align*}
	&\io (\chi+\rho^*(1-\chi))\big(V_R^k(p)+U^k[p]\big)(x,\tau) \dd x + (1+2k)\int_0^\tau\io \mu_R(p)|p|^{2k}|\nabla p|^2(x,t) \dd x\dd t \\
	& + \int_0^\tau\int_{\partial\Omega} \alpha(x)(p-p^*)\,p|p|^{2k}(x,t) \dd s(x)\dd t \\
	&\le C \int_0^\tau\io \left(|\dive u_t||p|^{1+2k}+(1+\hat\theta)|p|^{2k}+|p||p|^{1+2k}\right)(x,t) \dd x\dd t \\
	& + \io \left(V_R^k(p)+U^k[p] \right)\!(x,0) \dd x
\end{align*}
where, from Hypothesis \ref{h1}\,(vi),
\begin{align*}
	V_R^k(p)(x,t) &:= \int_0^{p(x,t)} f_R'(z)z|z|^{2k} \dd z \ge \int_0^p \frac{f^\flat}{(1+|z|)^{1+\nu}}\,z|z|^{2k} \dd z \ge \int_0^p \frac{f^\flat}{2\max\{1,|z|\}^{1+\nu}}\,z|z|^{2k} \dd z \\
	& =\frac{f^\flat}{2}\left(\int_0^1 z|z|^{2k} \dd z + \int_1^p z|z|^{2k-1-\nu} \dd z\right) \ge \frac{f^\flat}{1+2k-\nu} |p|^{1+2k-\nu} - C,\\
	U^k[p](x,t) &:= \int_0^\infty\int_0^{\play_r[p](x,t)} v|v|^{2k}\,\psi(r,v)\dd v\dd r \ge 0.
\end{align*}
Note also that, from Hypothesis \ref{h1}\,(vi) and an analogous version of \eqref{ini_play},
\begin{align*}
	&V_R^k(p)(x,0) = \int_0^{p^0(x)} f_R'(z)z|z|^{2k} \dd z \le \frac{f^\sharp}{2+2k}|p^0(x)|^{2+2k}, \\
	&U^k[p](x,0) = \int_0^\infty\int_0^{\play_r[p^0](x)} v|v|^{2k}\,\psi(r,v)\dd v\dd r \le C_\psi^* |p^0(x)|^{2+2k}.
\end{align*}
Moreover
$$|p|^{2k}|\nabla p|^2 = \frac{1}{(1+k)^2}|\nabla (p |p|^k)|^2.$$
Finally, by Young's inequality with conjugate exponents $\left(\frac{2+2k}{1+2k}\,,\,2+2k\right)$, we see that the boundary term is such that
\begin{align*}
	&\int_0^\tau\int_{\partial\Omega} \alpha(x)(p-p^*)\,p|p|^{2k} \dd s(x)\dd t \\
	&= \int_0^\tau\int_{\partial\Omega} \alpha(x)|p|^{2+2k} \dd s(x)\dd t - \int_0^\tau\int_{\partial\Omega} \alpha(x)p^*p|p|^{2k}(x,t) \dd s(x)\dd t \\
	&\ge 
	\frac{1}{2+2k} \int_0^\tau\int_{\partial\Omega} \alpha(x)|p|^{2+2k} \dd s(x)\dd t - \frac{1}{2+2k} \int_0^\tau\int_{\partial\Omega} \alpha(x)|p^*|^{2+2k} \dd s(x)\dd t.
\end{align*}
Hence, using also Hypothesis \ref{h1}\,(vii), we obtain
\begin{align*}
	&\frac{f^\flat}{1+2k-\nu} \io |p(x,\tau)|^{1+2k-\nu} \dd x + \mu^\flat\frac{1+2k}{(1+k)^2} \int_0^\tau\io |\nabla(p|p|^k)|^2 \dd x\dd t \\
	&\quad + \frac{1}{2+2k} \int_0^\tau\int_{\partial\Omega} \alpha(x)|p|^{2+2k} \dd s(x)\dd t \\
	&\le \frac{1}{2+2k} \int_0^\tau\int_{\partial\Omega} \alpha(x)|p^*|^{2+2k} \dd s(x)\dd t  + \left(\frac{f^\sharp}{2+2k}+ C_\psi^*\right) \io |p^0(x)|^{2+2k} \dd x  \\
	&\quad + C \int_0^\tau\io \left(|\dive u_t||p|^{1+2k}+(1+\hat\theta)\,(1+|p|)^{1+2k}+|p||p|^{1+2k}\right) \dd x\dd t.
\end{align*}
From Hypothesis \ref{h1}\,(iv) and (v) it follows that the above inequality is of the form
\begin{equation}\label{ee4}
	\begin{aligned}
		&\io |p(x,\tau)|^{1+2k-\nu} \dd x + \int_0^\tau\io |\nabla(p|p|^k)|^2 \dd x\dd t + \int_0^\tau\int_{\partial\Omega} \alpha(x)|p|^{2+2k} \dd s(x)\dd t \\
		&\qquad \le C(1+k)\left(C^k + \int_0^\tau\io |\tilde h||p|^{1+2k} \dd x\dd t\right)
	\end{aligned}
\end{equation}	
with a constant $C\ge 1$ independent of $R$ and $k$ and with
\be{thdef}
\tilde h = 1 + \hat\theta +|\dive u_t| +|p|.
\ee
We have $\tilde h \in L^2(\Omega\times(0,T))$ and
$$\|\tilde h\|_{L^2(\Omega\times(0,T))} \le C$$
independently of $R$ by virtue of the estimates \eqref{div} and \eqref{esti5} for $\dive u_t$, \eqref{esti_daf2} for $\hat\theta$ and \eqref{esti6} for $p$.
\smallskip

Since we deal with anisotropic spaces $L^q(0,T; L^r(\Omega))$, $q \ne r$, it is convenient to introduce for the norm of a function $v\in L^q(0,T; L^r(\Omega))$ the symbol
\begin{equation}\label{norqr}
	\|v\|_{r \ani q} := \left(\int_0^T |v(t)|_r^q \dd t\right)^{1/q}.
\end{equation}

For the function
\begin{equation}\label{ee4a}
	w_k(x,t) = p(x,t)|p(x,t)|^k
\end{equation}
we obtain from \eqref{ee4} using H\"older's inequality and Poincar\'e's inequality \eqref{poin} that
\begin{equation}\label{ee5a}
	\supess_{\tau \in (0,T)}|w_k(\tau)|_{s_k}^{s_k} + \int_0^\tau |w_k(t)|_{1;2}^2\dd t \le C(1+k)\left(C^k + \left(\int_0^T|w_k(t)|_{q_k}^{q_k}\dd t\right)^{1/2}\right)
\end{equation}
for all $\tau\in [0,T]$, with
$$s_k = \frac{1+2k-\nu}{1+k}\,, \qquad q_k = \frac{2+4k}{1+k}$$
and with a constant $C$ independent of $\tau$, $R$ and $k$. We now show that for a suitably chosen $k$, the right-hand side of \eqref{ee5a} is dominated by the left-hand side, which will imply a bound for the left-hand side. By the Gagliardo-Nirenberg inequality \eqref{gn} with $q=q_k$, $s=s_k$, $r=2$ and $N=3$ we have
\begin{equation}\label{gn2}
	|w_k(t)|_{q_k} \le C |w_k(t)|_{s_k}^{1-\delta_k} |w_k(t)|_{1;2}^{\delta_k}\,, \qquad \delta_k = \frac{\frac1{s_k} - \frac1{q_k}}{\frac1{s_k} -\frac16}.
\end{equation}
We now choose $k$ in such a way that $\delta_k q_k = 2$, that is, $3q_k = 6 + 2s_k$, which yields
$$k = 1-\nu, \quad s_k = \frac{3(1-\nu)}{2-\nu}, \quad q_k = \frac{6-4\nu}{2-\nu}, \quad q_k(1-\delta_k) = q_k-2 = \frac{2(1-\nu)}{2-\nu} = \frac23 s_k.$$
By Hypothesis \ref{h1} we have $s_k \ge 1$. Hence, by \eqref{gn2},
$$\int_0^T|w_k(t)|_{q_k}^{q_k}\dd t \le C \supess_{\tau \in (0,T)}|w_k(\tau)|_{s_k}^{(2/3)s_k} \int_0^T|w_k(t)|_{1;2}^2\dd t.$$
Since $k<1$, we conclude from \eqref{ee5a} that there exists a constant $C$ independent of $R$ such that, in particular,
$$
\supess_{\tau \in (0,T)} |w_k(\tau)|_{s_k} \le C\,,\quad \int_0^T|w_k(t)|_{q_k}^{q_k}\dd t \le C\,.$$
Invoking \eqref{ee4a}, we obtain for $p$ the estimates
\begin{equation}\label{ee8}
	\supess_{\tau \in (0,T)}|p(\tau)|_{3(1-\nu)} \le C\,, \quad \int_0^T|p(t)|_{6-4\nu}^{6-4\nu}\dd t \le C\,.
\end{equation}
We now distinguish two cases: $\nu \le 1/3$ and $\nu > 1/3$. For $\nu \le 1/3$ (that is, $3(1-\nu) \ge 2$) we have
\begin{equation}\label{ee8a}
	\supess_{\tau \in (0,T)}|p(\tau)|_{2} \le C\,.
\end{equation}
For $\nu > 1/3$ (that is, $3(1-\nu) < 2$) we use again the Gagliardo-Nirenberg inequality \eqref{gn} with $q=2$, $s=3(1-\nu)$, $r=2$ and $N=3$, obtaining
$$|p(t)|_2 \le C|p(t)|_{3(1-\nu)}^{1-\delta}|p(t)|_{1;2}^\delta\,, \qquad \delta = \frac{3\nu-1}{1+\nu}\,,$$
so that for
$$q_\nu = \frac{2(1+\nu)}{3\nu-1}$$
we have
$$\left(\int_0^T |p(t)|_2^{q_\nu} \dd t\right)^{1/q_\nu} \le C \supess_{t \in (0,T)}|p(t)|_{3(1-\nu)}^{1-\delta} \, \left(\int_0^T|p(t)|_{1;2}^2 \dd t\right)^{1/q_\nu}.$$
Hence by virtue of \eqref{esti6}, \eqref{norqr}, and \eqref{ee8} we obtain
\begin{equation}\label{eep}
	\|p\|_{2 \ani q_\nu} \le C
\end{equation}
with a constant $C>0$ independent of $R$, according to the notation \eqref{norqr}. Note that $q_\nu \ge 6$ thanks to Hypothesis \ref{h1} with $\nu \le 1/2$.


\subsection{Further estimates}\label{subsec_estu}

Now that we have obtained a suitable estimate for the capillary pressure in \eqref{eep}, we derive an analogous estimate for $\dive u_t$. To this aim we test \eqref{e2r} by $\psi = u_t$, which yields
$$
\begin{aligned}
&\io (\BB\nabla_s u_t:\nabla_s u_t)(x,t)\dd x \\
&\le 
\io \left(-P[\nabla_s u]:\nabla_s u_t + |p||\dive u_t| + \beta|\hat\theta-\theta_c||\dive u_t| +|g||u_t|\right)(x,t) \dd x.
\end{aligned}
$$
By \eqref{ivP}, Hypothesis \ref{h1}\,(v) and \eqref{esti5} we have
$$\io |P[\nabla_s u]|^2(x,t) \dd x \le C\left(1 + \int_0^T\io |\nabla_s u_t|^2(x,\tau) \dd x\dd\tau\right) \le C,$$
hence using Hypothesis \ref{h1}\,(i) and Young's inequality as in Subsection \ref{subsec_mechen} we conclude that the estimate
\begin{equation}\label{eu1}
	\io |\nabla_s u_t|^2(x,t)\dd x \le C\left(1 + \io \left(p^2 + \hat\theta^2\right)(x,t) \dd x\right)
\end{equation}
holds for a.\,e. $t\in(0,T)$ with a constant $C>0$ independent of $R$. We want to find and estimate for $\nabla_s u_t$ in the norm of $L^q(0,T; L^2(\Omega))$ for a suitable $q$. To this aim we apply the Gagliardo-Nirenberg inequality \eqref{gn} to $\hat\theta$ with the choices $q=r=2$, $s=1+b$ (with $b$ from Hypothesis \ref{h1}) and $N=3$. We obtain that, for $t\in(0,T)$,
$$|\hat\theta(t)|_2 \le C |\hat\theta(t)|_{1+b}^{1-\delta}|\hat\theta(t)|_{1;2}^{\delta}, \qquad \delta = \frac{3-3b}{5-b},$$
so that for
$$q_b = \frac{2(5-b)}{3-3b}$$
we have
$$\left(\int_0^T |\hat\theta(t)|_2^{q_b} \dd t\right)^{1/q_b} \le C \supess_{t \in (0,T)} |\hat\theta(t)|_{1+b}^{1-\delta} \left(\int_0^T |\hat\theta(t)|_{1;2}^2 \dd t\right)^{1/q_b}.$$
Hence by virtue of \eqref{esti_en3} and \eqref{esti_daf1} we get
\begin{equation}\label{eet}
	\|\hat\theta\|_{2\ani q_b} \le C
\end{equation}
independently of $R$. Note that our hypotheses on $b$ imply $q_b \ge 6$. Thus, coming back to \eqref{eu1} we have obtained that there exists $q := \min\{q_\nu,q_b\} \ge 6$ such that, thanks to \eqref{ee8a} or \eqref{eep} and \eqref{eet},
\begin{equation}\label{eeu}
	\|\nabla_s u_t\|_{2\ani q} \le C
\end{equation}
independently of $R$, according to the notation \eqref{norqr}.



We continue our analysis with the inequality \eqref{ee4} again. Unlike in Subsection \ref{subsec_estp}, we do not keep the exponent $k$ bounded, but we let $k \to \infty$ in a controlled way. As in \eqref{ee4a}, we define auxiliary functions $w_k = p|p|^k$ and rewrite \eqref{ee4} for $k \ge 1-\nu$ as
\begin{equation}\label{m10}
	\begin{aligned}
		&\io |w_k(x,\tau)|^{2a_k} \dd x + \int_0^\tau\io|\nabla w_k|^2 \dd x\dd t + \int_0^\tau\int_{\partial\Omega} \gamma(x) |w_k|^2 \dd s(x)\dd t\\
		&\qquad \le (1+k)\max\left\{C^{1+k},\int_0^\tau\io |\tilde h||w_k|^{2b_k} \dd x\dd t\right\}
	\end{aligned}
\end{equation}
with a constant $C\ge 1$ independent of $k$, with $\tilde h$ given by \eqref{thdef}, and with
$$a_k = \frac{1+2k-\nu}{2+2k}\,, \quad b_k = \frac{1+2k}{2+2k}\,.$$
It follows from \eqref{ee8a} or \eqref{eep}, \eqref{eet}, and \eqref{eeu} that $\tilde h \in L^6(0,T;L^2(\Omega))$ and
$$
	\|\tilde h\|_{2\ani 6} \le C.
$$
Repeating exactly the argument of the proof of \cite[Proposition 6.2]{gk} we obtain the following result.

\begin{propo}\label{prop_moser}
	Let Hypothesis \ref{h1} hold and let $(p,u,\theta,\chi)$ be a solution of \eqref{e1r}--\eqref{e4r} with the regularity from Proposition \ref{p1}. Then the function $p$ admits an $L^\infty$-bound independent of $R$, more precisely,
	\begin{equation}\label{mo1}
		|p(x,t)| \le C \left(({\bar\nu}H)^{\sigma/({\bar\nu}(\sigma - 1))} \sigma^{\sigma/({\bar\nu}(\sigma - 1)^2)}\right) =: R_\sigma
	\end{equation}
	for a.\,e. $(x,t) \in \Omega\times (0,T)$ with $\sigma = 19/18$, $H = \max\left\{1, \|\tilde h\|_{2\ani 6}\right\}$, and with positive constants $\bar \nu, C$ depending only on the data.
\end{propo}

The main consequence of Proposition \ref{prop_moser} is that, since we aim at taking the limit as $R \to \infty$ in \eqref{e1r}--\eqref{e4r}, we can restrict ourselves to parameter values $R > R_\sigma$, with $R_\sigma$ from \eqref{mo1}, so that the cut-off \eqref{fR}, \eqref{R}, \eqref{muR} is never active and $\gamma_R(p,\theta,\dive u) = \gamma(\hat\theta,\dive u)$. Hence we can rewrite \eqref{e1r}--\eqref{e4r} in the form
\begin{align}
	&
	\begin{aligned}
		&\int_{\Omega} \big((\chi+\rho^*(1-\chi))(f(p)+G_0[p]+\dive u)\big)_t\,\phi \dd x+\int_{\Omega} \frac{1}{\rho_W}\mu(p)\nabla p\cdot\nabla\phi \dd x \\
		&= \int_{\partial\Omega} \alpha(x)(p^*-p)\,\phi \dd s(x),
	\end{aligned}
	\label{e1moser}
	\\[4 mm]
	&\int_{\Omega} (P[\nabla_s u]+\BB\nabla_s u_t):\nabla_s\psi \dd x - \int_{\Omega} \big(p(\chi+\rho^*(1-\chi))+\beta(\hat\theta-\theta_c)\big)\,\dive\psi \dd x = \int_{\Omega} g\cdot\psi \dd x, \label{e2moser}
	\\[2 mm]
	&
	\begin{aligned}
		&\io \left(\mathcal{C}_V(\theta)_t\,\zeta+\kappa(\hat\theta)\nabla\theta\cdot\nabla\zeta\right) \dd x + \int_{\partial\Omega} \omega(x)(\theta-\theta^*)\,\zeta \dd s(x) \\
		& = \io\Bigg(\BB\nabla_s u_t:\nabla_s u_t+\|D_P[\nabla_s u]_t\|_*+\frac{1}{\rho_W}\mu(p)\,Q_R(|\nabla p|^2)+(\chi+\rho^*(1-\chi))|D_0[p]_t| \\
		& \ +\gamma(\hat\theta,\dive u)\chi_t^2-\left(\frac{L}{\theta_c}\chi_t+\beta\dive u_t\right)\hat\theta\Bigg)\zeta \dd x,
	\end{aligned}
	\label{e3moser}
\end{align}
\begin{align}
	&\gamma(\hat\theta,\dive u)\chi_t+\partial I_{[0,1]}(\chi) \ni (1-\rho^*)\left(\Phi(p)+p\,G_0[p]-U_0[p]+p\,\dive u\right)+L\left(\frac{\hat\theta}{\theta_c}-1\right) \ \text{a.\,e.} \label{e4moser}
\end{align}
for all test functions \(\phi, \zeta \in X\), \(\psi \in X_0\), with $\hat\theta = Q_R(\theta)$ and with initial conditions \eqref{ini}. In order to pass to the limit as $R \to \infty$, we still need to derive some higher order estimates.


\subsection{Higher order estimates for the capillary pressure}

Let us define
\begin{equation}\label{M}
	M(p) := \int_0^p \mu(z) \dd z
\end{equation}
for $p \in \real$, so that $\mu(p)\nabla p = \nabla M(p)$. We would like to test \eqref{e1moser} by $\phi = M(p)_t = \mu(p)p_t$ which, however, is not an admissible test function since $p_t \notin X$. Hence we choose a small $h>0$ and test by $\phi = \frac{1}{h}\left(M(p)(t)-M(p)(t-h)\right)$, where
$$\phi(x,t) = \frac{1}{h}\left(M(p)(t)-M(p)(t-h)\right)\!(x) := \frac{1}{h}\left(M(p(x,t))-M(p(x,t-h))\right),$$
with the intention to let $h \to 0$. We obtain
\begin{equation}\label{hop1}
	\begin{aligned}
		&\int_{\Omega} \big((\chi+\rho^*(1-\chi))(f(p)+G_0[p]+\dive u)\big)_t\,\frac{1}{h}\left(M(p)(t)-M(p)(t-h)\right) \dd x \\
		& + \int_{\Omega} \frac{1}{\rho_W}\nabla M(p)\cdot\nabla\left(\frac{1}{h}\left(M(p)(t)-M(p)(t-h)\right)\right) \dd x \\
		&= \int_{\partial\Omega} \alpha(x)(p^*-p)\,\frac{1}{h}\left(M(p)(t)-M(p)(t-h)\right) \dd s(x).
	\end{aligned}
\end{equation}
Concerning the second summand on the left-hand side of \eqref{hop1}, note that
\begin{align*}
	\nabla &M(p)(x,t)\cdot\nabla\left(\frac{1}{h}\left(M(p)(x,t)-M(p)(x,t-h)\right)\right) \\
	&\ge \frac{1}{2h} \left(|\nabla M(p)|^2(x,t)-|\nabla M(p)|^2(x,t-h)\right).
\end{align*}
We now deal with the boundary term. It holds
\begin{align*}
	&p^*(x,t)\,\frac{1}{h}\left(M(p)(x,t)-M(p)(x,t-h)\right) \\
	&= \frac{1}{h}\,\big(p^*(x,t)\,M(p)(x,t)-p^*(x,t-h)M(p)(x,t-h)\big) \\
	& \ - \frac{1}{h}\left(p^*(x,t)-p^*(x,t-h)\right)M(p)(x,t-h),
\end{align*}
To handle the term $p\,\frac{1}{h}\left(M(p)(t)-M(p)(t-h)\right)$, we use the inequality $F(y)-F(z) \le F'(y)(y-z)$ which holds for every convex function $F$ and every $y,z$. We interpret $M(p)(t)$ as $y$, $M(p)(t-h)$ as $z$ and $F'(y) = M^{-1}(M(p(t))) = M^{-1}(y)$. The function $M^{-1}$ is increasing, hence its antiderivative $F$ is convex. Thus
\begin{align*}
	p\left(M(p)(t)-M(p)(t-h)\right) \ge \int_{M(p)(t-h)}^{M(p)(t)} M^{-1}(z) \dd z = \int_{p(t-h)}^{p(t)} \xi\,M'({\xi}) \dd\xi.
\end{align*}
Defining
$$\hat\mu(p) := \int_0^p z\mu(z) \dd z$$
for $p \in \real$, we obtain
$$p(x,t)\,\frac{1}{h}\left(M(p)(x,t)-M(p)(x,t-h)\right) \ge \frac{1}{h}\left(\hat\mu(p)(x,t)-\hat\mu(p)(x,t-h)\right).$$
Thus \eqref{hop1} and the above estimates entail
\begin{align*}
	&\int_{\Omega} (\chi+\rho^*(1-\chi))(f(p)_t+G_0[p]_t)(x,t)\,\frac{1}{h}\left(M(p)(x,t)-M(p)(x,t-h)\right) \dd x \\
	& + \frac{1}{2\rho_W} \int_{\Omega} \frac{1}{h} \, \big(|\nabla M(p)|^2(x,t)-|\nabla M(p)|^2(x,t-h)\big) \dd x \\
	& + \int_{\partial\Omega} \alpha(x) \, \frac{1}{h} \left(\hat\mu(p)(x,t)-\hat\mu(p)(x,t-h)\right) \dd s(x) \\
	& - \int_{\partial\Omega} \alpha(x) \, \frac{1}{h} \, \big(p^*(x,t)\,M(p)(x,t)-p^*(x,t-h)M(p)(x,t-h)\big) \dd s(x) \\
	&\le - \int_{\Omega} (1-\rho^*)\chi_t\,(f(p)+G_0[p]+\dive u)(x,t)\,\frac{1}{h}\left(M(p)(x,t)-M(p)(x,t-h)\right) \dd x \\
	& - \int_{\Omega} (\chi+\rho^*(1-\chi))\,\dive u_t(x,t)\,\frac{1}{h}\left(M(p)(x,t)-M(p)(x,t-h)\right) \dd x \\
	& - \int_{\partial\Omega} \alpha(x) \, \frac{1}{h}\left(p^*(x,t)-p^*(x,t-h)\right)M(p)(x,t-h) \dd s(x).
\end{align*}
We are now ready to integrate in time from $h$ to some $\tau \in (0,T)$ and then let $h \to 0$. Note that estimate \eqref{esti4g} entails that the function $M(p)_t = \mu(p) p_t$ is in $L^2$, so that the convergence is strong in $L^2$. We obtain
\begin{align*}
	&\int_0^\tau \io (\chi+\rho^*(1-\chi))(f(p)_t+G_0[p]_t)\,\mu(p)p_t \dd x \dd t \\
	& + \frac{1}{2\rho_W} \io |\nabla M(p)|^2(x,\tau) \dd x - \frac{1}{2\rho_W} \io |\nabla M(p^0)|^2(x) \dd x \\
	& + \int_{\partial\Omega} \alpha(x) \left(\hat\mu(p)-p^*M(p)\right)(x,\tau) \dd s(x) - \int_{\partial\Omega} \alpha(x) \left(\hat\mu(p^0)(x)-p^*(x,0)M(p^0)(x)\right) \dd s(x) \\
	&\le - \int_0^\tau \int_{\Omega} (1-\rho^*)\chi_t\,(f(p)+G_0[p]+\dive u)\,\mu(p)p_t \dd x \dd t \\
	& - \int_0^\tau \int_{\Omega} (\chi+\rho^*(1-\chi))\,\dive u_t\,\mu(p)p_t \dd x \dd t - \int_0^\tau \int_{\partial\Omega} \alpha(x) \, p^*_tM(p) \dd s(x) \dd t.
\end{align*}
Combining \eqref{G0} with the identity \eqref{id_frak} for the play, we see that it holds $\mu(p)\,G_0[p]_t\,p_t \ge 0$, thus
$$(\chi+\rho^*(1-\chi))\left(f(p)_t+G_0[p]_t\right)\mu(p)p_t \ge \rho^*\mu^\flat\frac{f^\flat}{2\max\{1,R_\sigma\}^{1+\nu}}\,|p_t|^2$$
thanks to Hypothesis \ref{h1}\,(vi),\,(vii) and estimate \eqref{mo1}. Hence there exists a constant $c>0$ such that for every $t \in (0,T)$ we have
\begin{align*}
	&c \int_0^\tau \io |p_t|^2(x,t) \dd x \dd t + \frac{(\mu^\flat)^2}{2\rho_W}\io |\nabla p|^2(x,\tau) \dd x + \frac{\mu^\flat}{2}\int_{\partial\Omega} \alpha(x)p^2(x,\tau) \dd s(x) \\
	&\le C \Bigg( 1 + \int_0^\tau\io \left(|\chi_t|(1+|\dive u|)|p_t| + |\dive u_t||p_t|\right) \dd x \dd t + \int_0^\tau \int_{\partial\Omega} \alpha(x)|p^*_t||p| \dd s(x) \dd t \Bigg)
\end{align*}
thanks to Hypothesis \ref{h1}\,(v), (vi) and (vii), where we handled the boundary term on the left-hand side as in \eqref{eg7}. Arguing as for estimate \eqref{esti0g}, we obtain
\begin{equation}\label{esti13}
	|\chi_t(x,t)| \le \left|\frac{(1-\rho^*)\left(\Phi(p)+p\,G_0[p]-U_0[p]+p\,\dive u\right)+L\left(\hat\theta/\theta_c-1\right)}{\gamma(\hat\theta,\dive u)}\right| \le C
\end{equation}
for a.\,e. $(x,t) \in \Omega\times(0,T)$, this time independently of $R$ thanks to \eqref{mo1}. Thus, employing also Young's inequality, estimates \eqref{esti_en1}, \eqref{esti_en2}, \eqref{esti5} and Hypothesis \ref{h1}\,(iii) and (iv), we conclude that
\begin{align}
	\supess_{\tau \in (0,T)} \left(\io |\nabla p|^2(x,\tau) \dd x + \int_{\partial\Omega} \alpha(x)\,p^2(x,\tau)\dd s(x)\right) &\le C, \label{esti7} \\
	\int_0^T \io p_t^2(x,t) \dd x \dd t &\le C. \label{esti8}
\end{align}
By \eqref{esti_en1}, \eqref{esti5}, \eqref{esti6}, \eqref{esti8} and by comparison in equation \eqref{e1moser}, we see that the term $\Delta M(p)$ is bounded in $L^2(\Omega\times(0,T))$ independently of $R$. In terms of the new variable $\tilde{p}=M(p)$, the boundary condition in \eqref{boufull} is nonlinear, and from considerations similar to those used in the proof of \cite[Theorem 4.1]{kp} it follows
\begin{equation}\label{esti9}
	\|M(p)\|^2_{L^2(0,T;W^{2,2}(\Omega))} \le C.
\end{equation}
We thus may employ the Gagliardo-Nirenberg inequality \eqref{gn} with $s=r=2$, $N=3$ obtaining
$$|\nabla M(p)(t)|_q \le C \left(|\nabla M(p)(t)|_2 + |\nabla M(p)(t)|_2^{1-\delta}|\Delta M(p)(t)|_2^\delta\right), \qquad \delta = 3\left(\frac{1}{2}-\frac{1}{q}\right),$$
which holds for all $2<q<6$. Elevating to some power $s$ such that $\delta s=2$ and integrating in time yield
$$\int_0^T |\nabla p(t)|_q^s \dd t \le C \qquad \text{for } q \in (2,6) \quad \text{and} \quad \frac{1}{q} + \frac{2}{3s} = \frac{1}{2}\,,$$
thanks to estimates \eqref{esti7}, \eqref{esti9} and Hypothesis \ref{h1}\,(vii). In particular, for $s=4,\,q=3$ and $s=q=\frac{10}{3}$ we obtain, respectively,
\begin{equation}\label{esti11}
	\|\nabla p\|_{3 \ani 4} \le C, \qquad\qquad \|\nabla p\|_{10/3} \le C,
\end{equation}
according to the notation \eqref{norqr}.


\subsection{Higher order estimates for the displacement}

Let us consider equation \eqref{e2moser}. Setting
$$w(x,t) := p(\chi+\rho^*(1-\chi))+\beta(\hat\theta-\theta_c)(x,t) + \widehat{G}(x,t)$$
and arguing as for \eqref{nasut^p} we deduce
\begin{equation}\label{hod1}
	\begin{aligned}
		\io |\nabla_s u_t|^r(x,t) \dd x &\le C \io |\nabla_s u^0|^r(x) \dd x + Ct^{r-1} \int_0^t\io |\nabla_s u_t|^r(x,\tau) \dd x \dd\tau \\
		& + C \io |w|^r(x,t) \dd x \quad \text{a.\,e.}
	\end{aligned}
\end{equation}
Thus, by choosing $r=8/3+2b/3$ (with $b \in [1/2,1)$ from Hypothesis \ref{h1}) in the above inequality we obtain from \eqref{esti_daf2}, \eqref{mo1}, Hypothesis \ref{h1}\,(ii), (v) and Gr\"onwall's lemma
\begin{equation}\label{esti12}
	\|\nabla_s u_t\|_{8/3+2b/3} \le C.
\end{equation}
We now derive an estimate for $\nabla_s u_t$ in a suitable anisotropic Lebesgue space. To this aim we need to derive first an additional estimate for $\hat\theta$. We use Gagliardo-Nirenberg inequality \eqref{gn} with the choices $s=1+b$, $r=2$ and $N=3$ obtaining, for $t\in(0,T)$,
$$|\hat{\theta}(t)|_q \le C\left(|\hat\theta(t)|_{1+b} + |\hat\theta(t)|_{1+b}^{1-\delta}|\nabla\hat\theta(t)|_2^{\delta}\right), \qquad \delta = \left(\frac{1}{1+b}-\frac{1}{q}\right)\frac{6(1+b)}{5-b}\,,$$
which holds for all $1+b<q<6$. This yields, elevating to some power $s$ such that $\delta s = 2$ and integrating in time,
$$\int_0^T |\hat\theta(t)|_q^s \dd t \le C \qquad \text{for }q \in (1+b\,,\,6) \quad \text{and} \quad \frac{1}{q}+\frac{5-b}{3s(1+b)}=\frac{1}{1+b}$$
thanks to estimates \eqref{esti_en3}, \eqref{esti_daf1}. In particular, for $q=\frac{12(1+b)}{7+b}$ and $s=4$ we obtain
\begin{equation}\label{hod2}
	\|\hat\theta\|_{12(1+b)/(7+b) \ani 4} \le C.
\end{equation}
Note that
$$\frac{12(1+b)}{7+b} < \frac83 + \frac{2b}{3} \; \Leftrightarrow \; -7<b<2 \ \vee \ b>5,$$
which is certainly true under our hypotheses. Therefore, choosing $r=\frac{12(1+b)}{7+b}$ in \eqref{hod1} we obtain, thanks to \eqref{esti12} and Hypothesis \ref{h1}\,(v),
$$\io |\nabla_s u_t|^{12(1+b)/(7+b)}(x,t) \dd x \le C \left(1 + \io |w|^{12(1+b)/(7+b)}(x,t) \dd x\right)$$
for a.\,e. $t\in(0,T)$. Hypothesis \ref{h1}\,(ii) and estimates \eqref{mo1}, \eqref{hod2} then yield
\begin{equation}\label{esti14}
	\|\nabla_s u_t\|_{12(1+b)/(7+b) \ani 4} \le C
\end{equation}
independently of $R$.


\subsection{Higher order estimates for the temperature}\label{hot}

Note that \eqref{e2moser} with $\psi = u_t$ and \eqref{e4moser} entail, respectively,
\begin{align*}
	\int_{\Omega} \BB\nabla_s u_t:\nabla_s u_t \dd x &= - \io P[\nabla_s u]:\nabla_s u_t \dd x + \int_{\Omega} \big(p(\chi+\rho^*(1-\chi))+\beta(\hat\theta-\theta_c)\big)\,\dive u_t \dd x \\
	& \quad + \int_{\Omega} g \cdot u_t \dd x, \\
	\gamma(\hat\theta,\dive u)\chi_t^2 &= (1-\rho^*)\chi_t\,\Big(\Phi(p)+p\,G_0[p]-U_0[p]+p\,\dive u\Big)+L\chi_t\left(\frac{\hat\theta}{\theta_c}-1\right).
\end{align*}
Plugging these identities into \eqref{e3moser} we obtain
\begin{align}\label{hot1}
	&\io \left(\mathcal{C}_V(\theta)_t\,\zeta+\kappa(\hat\theta)\nabla\theta\cdot\nabla\zeta\right) \dd x + \int_{\partial\Omega} \omega(x)(\theta-\theta^*)\,\zeta \dd s(x) \nonumber\\
	& = \io\bigg(-P[\nabla_s u]:\nabla_s u_t+\big(p(\chi+\rho^*(1-\chi))-\beta\theta_c\big)\,\dive u_t+g \cdot u_t \nonumber\\
	& \ +\|D_P[\nabla_s u]_t\|_*+\frac{1}{\rho_W}\mu(p)\,Q_R(|\nabla p|^2)+(\chi+\rho^*(1-\chi))|D_0[p]_t| \nonumber\\
	& \ +\chi_t\,\Big((1-\rho^*)\,\big(\Phi(p)+p\,G_0[p]-U_0[p]+p\,\dive u\big)-L\Big)\bigg)\zeta \dd x \nonumber\\
	&=: \io \Gamma(x,t)\,\zeta \dd x
\end{align}
for every $\zeta \in X$, where $\Gamma(x,t)$ has the regularity of the worst term. Estimates \eqref{mo1}, \eqref{esti13}, \eqref{korn}, \eqref{ivP}, \eqref{enep2}, \eqref{eneg2} yield
$$|\Gamma| \le C(1+ |\nabla_s u|^2 + |\nabla_s u_t|^2 + |\nabla p|^2 + |p_t|),$$
which from \eqref{esti8}, \eqref{esti11}, \eqref{esti12}, \eqref{esti14} implies
\begin{equation}\label{hotA}
	\|\Gamma\|_{4/3+b/3} \le C, \qquad\qquad \|\Gamma\|_{6(1+b)/(7+b) \ani 2} \le C
\end{equation}
independently of $R$, with $b$ as in Hypothesis \ref{h1}.
\smallskip

Assume now that for some $p_0 \ge 8/3+2b/3$ we have proved
\begin{equation}\label{hothattheta}
	\|\hat\theta\|_{p_0} \le C.
\end{equation}
We know that this is true for $p_0 = 8/3+2b/3$ by virtue of \eqref{esti_daf2}. Set
\begin{equation}\label{r0}
	r_0 = \frac{1+b}{4+b}\,p_0,
\end{equation}
and set $\zeta = \hat\theta^{r_0}$ in \eqref{hot1}. We obtain
\begin{equation}\label{hot7}
	\io \left(\mathcal{C}_V(\theta)_t\,\hat\theta^{r_0}+\kappa(\hat\theta)\nabla\theta\cdot\nabla\hat\theta^{r_0}\right) \dd x + \int_{\partial\Omega} \omega(x)(\theta-\theta^*)\,\hat\theta^{r_0} \dd s(x) = \io \Gamma\,\hat\theta^{r_0} \dd x.
\end{equation}
It holds
$$\io \mathcal{C}_V(\theta)_t\,\hat\theta^{r_0} \dd x = \frac{\dd}{\dd t}\io F_{r_0}(\theta) \dd x$$
where
$$F_{r_0}(\theta) := \int_0^\theta c_V(s)(Q_R(s))^{r_0} \dd s.$$
Observe that
$$F_{r_0}(\theta) \ge \frac{\hat\theta^{r_0+1+b}}{r_0+1+b}$$
by Hypothesis \ref{h1}\,(viii). Moreover, Hypothesis \ref{h1}\,(ix) entails
$$\io \kappa(\hat{\theta})\nabla\theta\cdot\nabla\hat{\theta}^{r_0} \dd x = \io \kappa(\hat{\theta})\,r_0\,\hat{\theta}^{r_0-1}\nabla\theta\cdot\nabla\hat{\theta} \dd x \ge r_0\kappa^\flat \io \hat{\theta}^{r_0+a}|\nabla\hat{\theta}|^2 \dd x.$$
Concerning the boundary term, we use Young's inequality with exponents $\left(\frac{r_0+1}{r_0}\,,\,r_0+1\right)$, and obtain
\begin{align*}
	&\int_{\partial\Omega} \omega(x)(\theta-\theta^*)\hat{\theta}^{r_0} \dd s(x) \ge \int_{\partial\Omega}\omega(x)\hat{\theta}^{r_0+1} \dd s(x) - \int_{\partial\Omega}\omega(x)\theta^*\hat{\theta}^{r_0} \dd s(x) \\
	&\ge \int_{\partial\Omega} \omega(x)\hat{\theta}^{r_0+1} \dd s(x) - \frac{r_0}{r_0+1} \int_{\partial\Omega} \omega(x)\hat{\theta}^{r_0+1} \dd s(x) - \frac{1}{r_0+1} \int_{\partial\Omega} \omega(x)(\theta^*)^{r_0+1} \dd s(x) \\
	&\ge \frac{1}{r_0+1}\int_{\partial\Omega} \omega(x)\hat{\theta}^{r_0+1} \dd s(x) - C
\end{align*}
by Hypothesis \ref{h1}\,(iii) and (iv). We now integrate \eqref{hot7} in time $\int_0^\tau \dd t$ for some $\tau \in [0,T]$. Thanks to the choice \eqref{r0} and H\"older's inequality with exponents $\left(\frac{4+b}{3}\,,\,\frac{4+b}{1+b}\right)$, the right-hand side is such that
$$\int_0^\tau\io \Gamma\,\hat\theta^{r_0} \dd x\dd t = \int_0^\tau\io \Gamma\,(\hat\theta^{p_0})^{(1+b)/(4+b)} \dd x\dd t \le \|\Gamma\|_{(4+b)/3}\,\|\hat\theta\|_{p_0}^{r_0} \le C$$
by estimates \eqref{hotA}, \eqref{hothattheta}. Hence we have obtained
\begin{equation}\label{hot2}
	\begin{aligned}
		&\frac{1}{r_0+1+b}\io \hat\theta^{r_0+1+b}(x,\tau) \dd x + r_0 \int_0^\tau\io \hat{\theta}^{r_0+a}|\nabla\hat{\theta}|^2(x,t) \dd x\dd t \\
		& + \frac{1}{r_0+1}\int_0^\tau\int_{\partial\Omega} \omega(x)\hat{\theta}^{r_0+1}(x,t) \dd s(x)\dd t \le C.
	\end{aligned}
\end{equation}
We now denote
$$r = 1 + \frac{r_0+a}{2}\,, \quad s = \frac{r_0+1+b}{r}\,, \quad v = \hat\theta^{r}$$
and rewrite \eqref{hot2} as
\begin{equation}\label{hot3}
	\io |v|^s(x,\tau) \dd x + \int_0^\tau\io |\nabla v|^2(x,t) \dd x\dd t \le C(r_0+1+b).
\end{equation}
We now apply Gagliardo-Nirenberg inequality \eqref{gn} to $v(t)$, $t\in(0,T)$, with $r=2$ and $N=3$. Choosing $q$ in such a way that $\delta q = 2$, that is, $q = \frac{2}{3}s+2$, and integrating in time from $0$ to $T$ we obtain
\begin{align*}
	\|v\|_q &\le C\left(\sup_{t\in[0,T]}|v(t)|_s + \sup_{t\in[0,T]}|v(t)|_s^{(q-2)/q}\|\nabla v\|_2^{2/q}\right) \le C\left(\sup_{t\in[0,T]}|v(t)|_s + \|\nabla v\|_2\right).
\end{align*}
Estimate \eqref{hot3} yields
$$\sup_{t\in[0,T]}|v(t)|_s \le C(r_0+1+b)^{1/s}, \qquad \|\nabla v\|_2 \le C(r_0+1+b)^{1/2},$$
so that $\|v\|_q \le C(r_0+1+b)$. Coming back to the variable $\hat{\theta}$, we have proved that
$$\|\hat{\theta}\|_{p_1}\le C(r_0+1+b) \qquad \text{for }p_1 = rq = \frac{5(1+b)p_0}{3(4+b)}+\frac83+a+\frac{2b}{3}\,.$$
We now proceed by induction according to the rule
$$p_{j+1} = \frac{5(1+b)p_j}{3(4+b)}+\frac{8}{3}+a+\frac{2b}{3}, \qquad r_j = \frac{(1+b)p_j}{(4+b)}.$$
We have the implication
$$
p_j < \frac{(8+3a+2b)(4+b)}{7-2b} \ \Longrightarrow \ p_{j+1} > p_j.
$$
Hence, the sequence $\{p_j\}$ is increasing and
 $\lim_{j \to \infty} p_j = \frac{(8+3a+2b)(4+b)}{7-2b}$. It follows that choosing $\bar p = p_j$ for some $j$ sufficiently large we obtain
\be{rbar}
\bar{r} :=\frac{(1+b)\bar{p}}{(4+b)}>\hat{a},
\ee
with $\hat{a}$ from Hypothesis \ref{h1}\,(ix),
and using also \eqref{hot2} we obtain
\begin{equation}\label{hot4}
	\|\hat{\theta}\|_{\bar{p}} + \supess_{t \in (0,T)}|\hat\theta(t)|_{\bar{r}+1+b} \le C
\end{equation}
with $\bar p$ arbitrarily close to $\frac{(8+3a+2b)(4+b)}{7-2b}$. We now come back to \eqref{hot1}, which we test by $\zeta=\theta$ (note that this is an admissible choice by Proposition \ref{p1}). It holds
$$\io \mathcal{C}_V(\theta)_t\theta(x,t) \dd x = \io c_V(\theta)\theta\theta_t(x,t) \dd x = \frac{\dd}{\dd t} \io \left(\int_0^{\theta(x,t)} c_V(s)s \dd s\right) \dd x,$$
hence from Hypothesis \ref{h1}\,(ix) and \eqref{hotA} we obtain, after a time integration,
\begin{equation}\label{hot5}
	\begin{aligned}
		&\io \theta^{2+b}(x,\tau) \dd x + \int_0^T\io \kappa(\hat\theta)|\nabla\theta|^2(x,t) \dd x\dd t \\
		& + \int_0^T\int_{\partial\Omega} \omega(x)\theta^2(x,t) \dd s(x)\dd t \le C\left(1+\|\theta\|_{(4+b)/(1+b)}\right).
	\end{aligned}
\end{equation}
Using the Gagliardo-Nirenberg inequality \eqref{gn} again with $q=\frac{4+b}{1+b}$, $s=1+b$ (note that $1+b<\frac{4+b}{1+b}<6$ under our hypotheses), $r=2$ and $N=3$ we have that, for each fixed $t\in(0,T)$,
$$|\theta(t)|_{(4+b)/(1+b)} \le C\left(1 + |\nabla\theta(t)|_2^\delta\right)$$
with $\delta=\frac{6(3-b^2-b)}{(5-b)(4+b)}$ and where we used estimate \eqref{esti_en3}. Raising to the power $(4+b)/(1+b)$ and integrating $\int_0^T \dd t$ we get
$$\|\theta\|_{(4+b)/(1+b)} \le C\left(1 + \left(\int_0^T\io |\nabla\theta|^2 \dd x\dd t\right)^{\delta/2}\right) \le C\left(1 + \left(\int_0^T\io \kappa(\hat\theta)|\nabla\theta|^2 \dd x\dd t\right)^{\delta/2}\right).$$
Plugging this back into \eqref{hot5} and using Young's inequality we deduce
\begin{equation}\label{ho7}
	\io \theta^{2+b}(x,\tau) \dd x + \int_0^T\io \kappa(\hat\theta)|\nabla\theta|^2(x,t) \dd x\dd t + \int_0^T\int_{\partial\Omega} \omega(x)\theta^2(x,t) \dd s(x)\dd t \le C.
\end{equation}
This enables us to derive an upper bound for the integral $\int_\Omega \kappa(\hat{\theta})\nabla\theta\cdot\nabla\zeta \dd x$, which we need for getting an estimate for $\theta_t$ from equation \eqref{hot1}. By H{\"o}lder's inequality and Hypothesis \ref{h1}\,(ix) we have that
\begin{align}
	\int_\Omega |\kappa(\hat{\theta})\nabla\theta\cdot\nabla\zeta| \dd x &= \int_\Omega |\kappa^{1/2}(\hat{\theta})\nabla\theta\cdot\kappa^{1/2}(\hat{\theta})\nabla\zeta| \dd x \nonumber\\
	&\le C\left(\int_\Omega \kappa(\hat{\theta})|\nabla\theta|^2 \dd x\right)^{1/2} \left(\int_\Omega \max\{1,\hat{\theta}^{1+\hat{a}}\}|\nabla\zeta|^2 \dd x\right)^{1/2}. \label{hot6}
\end{align}
Let us now choose $\hat{q}>1$ such that $(1+\hat{a})\hat{q} = 1+\bar{r}+b$, where $\bar{r}$ is defined in \eqref{rbar}. Note that such a $\hat{q}$ exists since $1+\bar{r}+b > 1+\hat{a}+b > 1+\hat{a}$. Defining
\begin{equation}\label{qstar}
	q^* := \frac{2\hat{q}}{\hat{q}-1} = 2 + \frac{2}{\hat{q}-1} > 2\,,
\end{equation}
we get from H{\"o}lder's inequality with conjugate exponents $\left(\hat{q}\,,\,\frac{q^*}{2}\right)$ that
$$\int_\Omega \hat{\theta}^{1+\hat{a}}|\nabla\zeta|^2 \dd x \le \left(\int_\Omega \hat{\theta}^{1+\bar{r}+b} \dd x\right)^{1/\hat{q}} \left(\int_\Omega |\nabla\zeta|^{q^*} \dd x\right)^{2/q^*} \le C\left(\int_\Omega |\nabla\zeta|^{q^*} \dd x\right)^{2/q^*}$$
by virtue of \eqref{hot4}. Inequality \eqref{hot6} then yields the bound
$$\int_\Omega |\kappa(\hat{\theta})\nabla\theta\cdot\nabla\zeta| \dd x \le C\left(\int_\Omega \kappa(\hat{\theta})|\nabla\theta|^2 \dd x\right)^{1/2} \left(\int_\Omega |\nabla\zeta|^{q^*} \dd x\right)^{1/q^*}.$$
Hence, by \eqref{ho7},
$$\int_0^T\int_\Omega |\kappa(\hat{\theta})\nabla\theta\cdot\nabla\zeta| \dd x\dd t \le C\|\zeta\|_{L^2(0,T;W^{1,q^*}(\Omega))}.$$
From \eqref{hotA} it follows that testing with $\zeta \in L^2(0,T;W^{1,q^*}(\Omega))$ is admissible, in the sense that the term $\Gamma\zeta$ is integrable. This is obvious if $q^* \ge 3$. For $q^*<3$ the space $W^{1,q^*}(\Omega)$ is embedded in $L^{q^*_S}(\Omega)$ with
$$\frac1{q^*_S} = \frac1{q^*} - \frac13  = \frac16 - \frac1{2\hat q} ,$$
so that $\frac3{4+b} + \frac1{q^*_S} \le 1$.  We thus obtain from \eqref{hot1} that
\begin{equation}\label{esti15}
	\int_0^T\int_\Omega \theta_t\zeta \dd x\dd t \le C \|\zeta\|_{L^2(0,T;W^{1,q^*}(\Omega))}.
\end{equation}


\section{Passage to the limit as $R \to \infty$}\label{sec_lim}

In this section we conclude the proof of Theorem \ref{t1} by passing to the limit in \eqref{e1moser}--\eqref{e4moser} as $R \to \infty$. Most of the convergences can be handled as at the end of Subsection \ref{subsec_estn}, hence we focus here on the main differences.

Let $R_i \nearrow \infty$ be a sequence such that $R_i > R_\sigma$, with $R_\sigma$ as in \eqref{mo1}, and let $(p,u,\chi,\theta) = (p^{(i)},u^{(i)},\chi^{(i)},\theta^{(i)})$ be solutions of \eqref{e1moser}--\eqref{e4moser} corresponding to $R = R_i$, with $\hat\theta = \hat\theta^{(i)} = Q_{R_i}(\theta^{(i)})$ and test functions $\phi,\zeta \in X$, $\psi \in X_0$. Our aim is to check that at least a subsequence converges as $i \to \infty$ to a solution of \eqref{e1w}--\eqref{e4w} with test functions $\phi \in X$, $\psi \in X_0$ and $\zeta \in X_{q^*}$.

\medskip

First, for the capillary pressure $p = p^{(i)}$ we have the estimates \eqref{mo1}, \eqref{esti7}, \eqref{esti8}, \eqref{esti9} and \eqref{esti11}, which imply that, passing to a subsequence if necessary,
$$
\begin{array}{rcll}
	p^{(i)}_t &\to& p_t &\quad \text{weakly in }L^2(\Omega \times (0,T)), \\
	p^{(i)} &\to& p &\quad \text{strongly in }L^q(\Omega;C[0,T]) \ \text{for all }q \in [1,\infty), \\
	\nabla p^{(i)} &\to& \nabla p &\quad \text{strongly in }L^q(\Omega\times(0,T);\real^3) \ \text{for all }q \in \left[1,\frac{10}{3}\right),
\end{array}
$$
by compact embedding. We easily show that
\begin{equation}\label{conv1}
	Q_{R_i}(|\nabla p^{(i)}|^2) \to |\nabla p|^2 \quad \text{strongly in }L^q(\Omega\times(0,T);\real^3) \ \text{for all }q \in \left[1\,,\,\frac{5}{3}\right).
\end{equation}
Indeed, let $\Omega_T^{(i)} \subset \Omega\times(0,T)$ be the set of all $(x,t) \in \Omega\times(0,T)$ such that $|\nabla p^{(i)}(x,t)|^2 > R_i$. By \eqref{esti11} we have
$$C \ge \int_0^T\io |\nabla p^{(i)}(x,t)|^{10/3} \dd x\dd t \ge \iint_{\Omega_T^{(i)}} |\nabla p^{(i)}(x,t)|^{10/3} \dd x\dd t \ge |\Omega_T^{(i)}| R_i^{5/3},$$
hence $|\Omega_T^{(i)}| \le C R_i^{-5/3}$. For $q<\frac53$ we use H\"older's inequality to get the estimate
\begin{align*}
	\int_0^T\io &\left|Q_{R_i}(|\nabla p^{(i)}|^2)-|\nabla p^{(i)}|^2\right|^q \dd x\dd t = \iint_{\Omega_T^{(i)}} \left|R_i-|\nabla p^{(i)}|^2\right|^q \dd x\dd t \le \iint_{\Omega_T^{(i)}} |\nabla p^{(i)}|^{2q} \dd x\dd t \\
	&\le \left(\iint_{\Omega_T^{(i)}} |\nabla p^{(i)}|^{10/3} \dd x\dd t\right)^{3q/5} |\Omega_T^{(i)}|^{1-3q/5} \le C R_i^{-(5-3q)/3},
\end{align*}
and \eqref{conv1} follows.

\medskip

For the temperature $\theta = \theta^{(i)}$ we proceed in a similar way. By estimates \eqref{ho7} and \eqref{esti15} we obtain
$$
\begin{array}{rcll}
	\nabla\theta^{(i)} &\to& \nabla\theta &\quad \text{weakly in }L^2(\Omega\times(0,T);\real^3), \\
	\theta^{(i)}_t &\to& \theta_t &\quad \text{weakly in }L^2(0,T;W^{-1,q^*}(\Omega)), \\
	\theta^{(i)} &\to& \theta &\quad \text{strongly in }L^2(\Omega\times(0,T)),
\end{array}
$$
where for the last convergence we exploited \cite[Theorem 5.1]{lions} and the embedding $W^{-1,q^*}(\Omega) \hookrightarrow W^{-1,2}(\Omega)$ (recall that $q^*>2$). Furthermore, estimate \eqref{hot4} entails that $\hat\theta^{(i)}$ are uniformly bounded in $L^q(\Omega\times(0,T))$ for every $q<\frac{(8+3a+2b)(4+b)}{7-2b}$. Hence a similar argument as above yields that
$$\hat\theta^{(i)} = Q_{R_i}(\theta^{(i)}) \to \theta \quad \text{strongly in }L^q(\Omega\times(0,T)) \ \text{for all }q \in \left[1\,,\,\frac{(8+3a+2b)(4+b)}{7-2b}\right).$$
Estimate \eqref{esti9} and the Sobolev embeddings yield an inequality similar to \eqref{diff}, but with a constant independent of both $\eta$ and $R$. This is enough to obtain that $\nabla_s u^{(i)} \to \nabla_s u$, $\nabla_s u^{(i)}_t \to \nabla_s u_t$ strongly in $L^2(\Omega;C([0,T];\tens))$ and in $L^2(\Omega\times(0,T);\tens)$, respectively. The strong convergences $\chi^{(i)} \to \chi$, $\chi^{(i)}_t \to \chi_t$ then follow as at the end of Subsection \ref{subsec_gal}, as well as the convergence of the hysteresis terms.

Therefore the limit as $i \to \infty$ yields a solution to \eqref{e1w}--\eqref{e4w}, and the proof of Theorem \ref{t1} is completed.

%

\end{document}